\documentclass[11pt]{amsart}
\usepackage{amsmath,amssymb, graphics, amscd,latexsym}
\makeatletter

\newtheorem{Theorem}{Theorem}
\newtheorem{Lemma}[Theorem]{Lemma}
\newtheorem{Corollary}[Theorem]{Corollary}
\newtheorem{Proposition}[Theorem]{Proposition}

\newtheorem{Remark}[Theorem]{Remark}

\newcommand{\eps}{\varepsilon}

\newcommand\Si{\Sigma}

\newcommand\De{\Delta}

\newcommand\cM{\mathcal  M}

\newcommand\cC{\mathcal  C}

\newcommand\ord{\text{ord}\/}

\newcommand\modulo{\text{modulo}\/}

\newcommand\Coeff{\text{Coeff}\/}

\newcommand\degree{\text{degree}\/}

\begin{document}
\title[Zariski pairs on sextics II
]
{  Zariski pairs on sextics II}

\author
[M. Oka ]
{Mutsuo Oka }
\address{\vtop{
\hbox{Department of Mathematics}
\hbox{Tokyo  University of Science}
\hbox{1-3 Kagurazaka,Shinjuku-ku}
\hbox{Tokyo 162-8601}
\hbox{\rm{E-mail}: {\rm oka@ma.kagu.tus.ac.jp}}
}}
\keywords{ torus type, Zariski pair, flex points, conical flex points}

\maketitle
\begin{center}
{Dedicated to Professor Kyoji Saito for his 60th birthday}
\end{center}
\pagestyle{headings}

\vspace{1cm}
\section{Introduction}
We continue to study Zariski pairs in  sextics.
In this paper, we study Zariski pairs of sextics which are not
irreducible.
 The idea of the
construction
of Zariski partner sextic for reducible cases is quit different from the
irreducible case. It is crucial to take the geometry of the components
and their mutual intersection data into account. When there is a line
component,
flex geometry (i.e., linear geometry) is concerned to the geometry of sextics
of torus type and non-torus type.
When there is no linear components, the geometry
is more difficult to distinguish  sextics of torus type. For this
reason, we introduce the notion of conical flexes.

We have observed  in \cite{OkaAtlas} that the case $\rho(C,5)=6$ is  critical in the sense that
the Alexander polynomial $\De_C(t)$ can be either trivial or non-trivial
for sextics. If $\rho(C,5)>6 $ (resp. $\rho(C,5)<6$), the Alexander
polynomial
is not trivial (resp. trivial) (\cite{OkaAtlas}).  For the definition of
$\rho(C,5)$-invariant, see \cite{OkaAtlas}.
Thus we concentrate ourselves in this paper the case $\rho(C,5)=6$.
In \cite{Reduced}, we have classified the possible
configurations
for reducible sextics of torus type. In particular, the configurations
with $\rho(C,5)=6$ are given as in Theorem 1 below.
Hereafter we use the same notations as \cite{OkaAtlas} for
 denoting component types. For example, $C=B_1+B_5$ implies that $C$ has 
a linear component $B_1$ and a quintic component $B_5$.
 We denote the configuration
of the singularities of $C$ by $\Si(C)$.
\begin{Theorem}{\rm(\cite{Reduced})}\label{ListTorus}
Assume that $C$ is a   reducible sextic of torus type  with
 $\rho(C,5)=6$
and only simple singularities.
Let $\Si_{in}$ be the inner singularities. Then the possible
  configurations of simple singularities
are as follows.

  \begin{enumerate}
  \item  $\Si_{in}=[A_5,4A_2]:$  $C=B_5+B_1$ and 
 $\Si(C)=[A_5,  4A_2, 2A_1],\,  [A_5, 4A_2, 3A_1]$,
$[A_5, 4A_2, 4A_1]$.
  \item $\Si_{in}=[2A_5,2A_2]${\rm :}
\begin{enumerate}
\item $C=B_1+B_5${\rm :}
 $\Si(C)=[2A_5,2A_2,2A_1],  [2A_5,2A_2,3A_1]$.
\item $C= B_1+B_1'+\, B_4${\rm :}
$\Si(C)=[2A_5,2A_2,3A_1],\,  [2A_5,2A_2,4A_1]$.
\item $C= B_2+B_4${\rm :}
 $ \Si(C)=[2A_5,2A_2,2A_1]$, $[2A_5,2A_2,3A_1]$.
\item $C=B_3+B_3'${\rm :}
$ \Si(C)=[2A_5,2A_2,3A_1]$.
\end{enumerate}
  \item $\Si_{in}=[E_6,A_5,2A_2]${\rm :} $C=B_1+B_5$,
$ \Si(C)=[E_6, A_5, 2A_2, 2A_1]$, $[E_6, A_5, 2A_2, 3A_1]$.
  \item  $\Si_{in}=[3A_5]${\rm :}
\begin{enumerate}
\item $C=B_1+B_5${\rm :}
  $\Si(C)=[3A_5,  2A_1]$.
\item $C=B_2+B_4${\rm :} $\Si(C)=[3A_5,2A_1]$.
\item $C=B_1+B_1'+B_4${\rm :}
  $\Si(C)=[3A_5, 3A_1]$.
\item $C=B_3+B_3'${\rm :}  $\Si(C)=[3A_5]$, $[3A_5, A_1]$,
 $[3A_5, 2A_1]$.
\item $C=B_1+B_2+B_3${\rm :}
$ \Si(C)=[3 A_5, 2A_1],  [3 A_5,  3A_1]$.
\item $C=B_1+B_1'+B_1''+B_3${\rm :}  $\Si(C)=[3A_5,3A_1]$, $[3A_5,4A_1]$,
\item $C=B_2+B_2'+B_2''${\rm :} $\Si(C)=[3A_5,3A_1]$.
\end{enumerate}
  \item $\Si_{in}=[2A_5,E_6]${\rm :}
\begin{enumerate}
\item $C=B_1+B_5${\rm :}
$  \Si(C)=[E_6,2A_5,2A_1]$.
\item $C=B_2+B_4${\rm :}
  $\Si(C)=[E_6,2A_5,2A_1]$.
\item $C=B_1+B_1'+B_4${\rm :}
$\Si(C)=[E_6,2A_5,3A_1]$.
\end{enumerate}
  \item $\Si_{in}=[A_8,A_5,A_2]${\rm :}  $C=B_1+B_5$,
$  \Si(C)=[A_8, A_5,A_2, 2A_1]$, $[A_8, A_5,A_2,  3A_1]$.
  \item $\Si_{in}=[A_{11},2A_2]${\rm :}
\begin{enumerate}
\item $C=B_2+B_4${\rm :}
 $\Si(C)=[A_{11},2A_2,2A_1]$, $[A_{11},2A_2, 3A_1]$.
\item $C=B_3+B_3'${\rm :}
$\Si(C)=[A_{11},2A_2,3A_1]$.
\end{enumerate}
  \item $\Si_{in}=[A_{11},A_5]${\rm :}
\begin{enumerate}
\item  $C=B_1+B_5${\rm :}
$\Si(C)=[A_{11},A_5, 2A_1]$.
\item $C=B_2+B_4${\rm :}
$\Si(C)=[A_{11},A_5,2A_1]$.
\item $C=B_3+B_3'${\rm :}
$\Si(C)=[A_{11},A_5]$, $[A_{11},A_5,A_1]$,
$[A_{11},A_5,2A_1]$,
\item $C=B_1+B_2+B_3${\rm :}
$\Si(C)=[A_{11},A_5,2A_1]$, $[A_{11},A_5,3A_1]$,
\end{enumerate}
  \item $\Si_{in}=[A_{17}]${\rm :}  $C=B_3+B_3'$,
$\Si(C)=[A_{17}]$, $[A_{17},A_1]$,  
$[A_{17},2A_1]$.
\end{enumerate}

\end{Theorem}

Our main result  in this paper is:
\begin{Theorem}\label{Main result}
There are Zariski partner sextics with the above configurations
with the following exceptions:
\begin{enumerate}
\item $\Si(C)=[A_5,4A_2,4A_1]$ with $C=B_5+B_1$.
\item $\Si(C)=[2A_5,2A_2,4A_1]$ with $C=B_1+B_1'+B_4$.
\item $\Si(C)=[E_6,A_5,2A_2,3A_1]$ with $C=B_5+B_1$.
\item $\Si(C)=[3A_5,4A_1]$ with $C=B_3+B_1+B_1'+B_1''$.
\item $\Si(C)=[E_6,2A_5,3A_1]$ with $C=B_4+B_1 +B_1'$.
\end{enumerate}
\end{Theorem}
The non-existence  of sextics of non-torus type
with  the above exceptional configurations will be explained
by flex  geometry.
The existence will be also explained by the flex geometry
for those which has a line components and by conical flex geometry
for the component type $B_4+B_2,\, B_2+B_2'+B_2''$.
\begin{Remark}
(1) In the list of Theorem \ref{Main result}, there are certainly several
 cases
which are already known. For example, the configuration $C=B_3+B_3' $
with one singularity $A_{17} $ is given by Artal \cite{Artal}.

\noindent
(2) In this paper, we only studied possible Zariski pairs of reducible
 sextics
$(C,C')$ where $C$ is of torus type and $C'$ is not of torus type.
On the other hand, the possibility of Zariski pairs among reducible
 sextics
of the same class is not discussed here. Several examples are known 
among reducible sextics of non-torus type. For such cases, Alexander 
polynomials can not distinguish the differnece. See 
papers  \cite{A-R-C-T, AT1,AT2,ARCI}

\end{Remark}
 
\section{Reducible sextics of non-torus type}

To compute explicit polynomials defining reducible sextics, it is not 
usually easy to look for  special degenerations into several irreducible components
starting from the generic sextics
$\sum_{i+j\le 6} a_{ij}x^i\,y^j$.
Recall that we have classified all possible reducible simple configurations
in \cite{Reduced} and it is easier to start from a fixed reducible
decomposition. In fact, the geometry of the configuration of a
reducible sextic
depends
very much on the geometry of each components.
A smooth point $P\in C$ is called a {\em flex} point
if the intersection multiplicity of the tangent line and $C$
at $P$ is strictly greater than 2.
First we recall the following fact for flex points (\cite{NambaBook,Okadual}).
\begin{Lemma}\label{flex-number} Let $C: F(X,Y,Z)=0$ be an irreducible
 plane curve
of degree $n$
 with singularities $\{P_1,\dots, P_k\}$. Then the number of flexes
 $\iota(C)$
is given by
\[
 \iota(C)=3\,n\,(n-2)\,-\,\sum_{i=1}^k\,\eps(P_i;C)
\]
where the second term $\eps(P_i;C)$ is the flex defect and given by the
 local intersection number of $C$ and the hessian curve of $C$ at $P_i$.
\end{Lemma}
Generic flex defect of simple singularities we use are
\begin{eqnarray}\label{flex-formula}
 \eps(A_1)=6,\,\eps(A_2)=8,\, \eps(A_{3\iota-1})=9\iota,\, (\iota\ge 2),
\, \eps(E_6)=22
\end{eqnarray}

Recall that flex  points of a curve are described by the hessian 
of the defining homogeneous equation. When we have an affine
 equation $C:\, f(x,y)=0$, flex points  in $\bf C^2$  are described by 
$f(x,y)=flex_f(x,y)=0$ (\cite{Okadual})
where \[
       flex_f(x,y)\, :=\, f_{xx}\,
       f_y^2-2\,f_{xy}\,f_x\,f_y+ \,f_{yy}\,f_x^2
      \]
This is an easy way to check flex points from the affine equation.

A sextic $C$ is  of (2,3)-torus type
if we can take a defining polynomial of the form
 $f_2(x,y)^3+f_3(x,y)^2=0$ where $\degree\, f_j=j$.
The intersections $f_2=f_3=0$ are singular points of $C$
and we call them {\em inner singularities}.
For a given sextic $C$ of torus type whose  singularities are
simple, the possible inner singularities are

$(\sharp):\quad \{A_2,A_5,A_8,A_{11},A_{14},A_{17},E_6\}$.
\newline
A convenient criterion  for $C$ to be  of torus
type 
is the existence
a certain conic $C_2$ such that 
$C_2\cap C\subset \Si(C)$ (Tokunaga's criterion \cite{Tokunaga-torus},
Lemma 3, \cite{ZariskiPairsI}).

A sextic of torus type $C$
is called {\em of linear torus type} if  the conic polynomial $f_2$ can be written
as  $f_2(x,y)=\ell(x,y)^2$  for some linear form $\ell(x,y)$ (\cite{OkaAtlas}). 
A  sextic of linear torus type can have only
$A_5,\, A_{11},\, A_{17}$ as inner singularties and 
the location of these singularities are  colinear.

The proof of Theorem \ref{Main result} is done by giving explicit
examples. 
For the better understanding of the existence or non-existence
of the Zariski pairs, we divide the above
configurations into the following classes.
\begin{enumerate}
\item $C$ has a quintic component. The corresponding component type is
      $B_5+B_1$ and $B_1$ is a flex tangent line. 
\item $C$ has a quartic component. There are two subcases.
\begin{enumerate}
\item $C=B_4+B_1+B_1'$. In this case,
two line components are flex tangent lines.
 \item $C=B_4+B_2$.
\end{enumerate}
\item $C$ has a cubic component. There are two subcases.
\begin{enumerate}
\item  Sextics of linear torus type.
\item  Sextics, not of linear torus type.
\end{enumerate}
\item $C=B_2+B_2'+B_2''$.
\end{enumerate}

\section{Configuration coming from quintic flex geometry}
Let $B_5$ be an irreducible quintic 
and let $P$ be a flex point of
$B_5$. 
We  denote the tangent line at $P$ by $L_P$.
We say that $P$ is   {\em a flex  of torus type} (respectively {\em a flex of 
non-torus type}) if $B_5\cup L_P$ is a sextic of torus type (resp.
of non-torus type).
 The following configurations are  mainly  related to the flex geometry of
 certain quintics. (By 'flex geometry', we mean the geometry of the tangent lines at the flex points and 
the curve.)
Recall that $\Si(B_5)$ is the configuration of the
 singularities of $B_5$.  Let $\iota$ be the number of flex points on $B_5$.

\begin{enumerate}
\item $C=B_5+B_1$ with $\Si(C)=[A_5,4A_2,kA_1],\,k=2,\,3,\, 4$. Then
      $\Si(B_5)=[4A_2,(k-2)A_1]$ for $k=2,\,3,\, 4$ and
$\iota= 13,\, 7,\, 1$  respectively. 
\item $C=B_5+B_1$ with $\Si(C)=[2A_5,2A_2,kA_1],\,k=2,3$.
Then  $\Si(B_5)=[A_5,2A_2,(k-2)A_1]$ and $\iota=11,\, 5$ respectively.
\item $C=B_5+B_1$ with $\Si(C)=[E_6,A_5,2A_2,kA_1],\, k=2,\,3$.
Then
$\Si(B_5)=[E_6,2A_2,(k-2)A_1]$ and $\iota=7,\, 1$ respectively
      for $k=2,3$. The case $k=3$ corresponds to sextics of torus type.
\item $C=B_5+B_1$ with $\Si(C)=[E_6,2A_5,2A_1]$.
The quintic $B_5$ has $\Si(C)=[E_6,A_5]$ and $\iota=5$.
\item $C=B_5+B_1$ with $\Si(C)=[3A_5,2A_1]$. 
Then $\Si(B_5)=[2A_5]$ and $\iota=9$.
\item $C=B_5+B_1$ with $\Si(C)=[A_8,A_5,A_2,kA_1],\,k=2,3$.
Then  $\Si(B_5)=[A_8,A_2,(k-2)A_1]$ and $\iota=10,\,4$ for $k=2,3$.
\item $C=B_5+B_1$ with $\Si(C)=[A_{11},A_5,2A_1]$.
Then $\Si(B_5)=[A_{11}]$ and $\iota=9$.
\end{enumerate}
We are going to show the stronger assertion for the above
 configurations:
{\em the Zariski partner sextic of non-torus type
are simply given by replacing  the flex line components $B_1$ for the
above cases, if $B_5$ has at least two flex points.}

Let $\Xi$ be a configuration of singularities on $B_5$, which is one of
the
 above list. Let $\cM(\Xi;5)$  be {\em the configuration space} of quintics
 $B_5$  such that $\Si(B_5)=\Xi$.
We considere it as a topological subspace of the space of quintics.
 For our purpose,
it is enough to consider the marked configuration subspace $\cM(\Xi;5)'$
which consists of the pair $(B_5,P)$, where    $B_5\in \cM(\Xi;5)$
 and $P$ is a  flex point of torus type. 
The following describes the existence of sextics of
non-torus type with the above configurations.
\begin{Theorem}\label{quintic-flex} Let $\Xi$ be a configuration of singularities on $B_5$, which is one of
the
 above list.
1. The configuration  subspace  $\cM(\Xi;5)'$ is connected for each $\Xi$.
\newline
2. For each  $\Xi\ne[4A_2,2A_1],\,[ E_6,2A_2,A_1]$ and $B_5\in \cM(\Xi;5)'$,
a Zariski
 pair sextics
are given as $\{B_5\cup L_P,\,B_5\cup L_Q\}$
where $P$ and $Q$ are
 flex points of torus-type and  of non-torus type respectively.
\newline 
3. For these two exceptional cases,
we have the equality
$\cM(\Xi;5)'=\cM(\Xi;5)$ and  a quintic $B_5\in\cM(\Xi;5) $ does not  contain any
flexes of non-torus type. \newline
\end{Theorem}

\begin{Remark}
Let $\iota_t,\,\iota_{nt}$  be the respective number of flex
 points of torus type
and of non-torus type on a generic $B_5\in \cM(\Xi;5)'$.
We do not need the precise number $\iota_t,\, \iota_{nt}$ for our
 purpose. The sum $\iota=\iota_t+\iota_{nt}$ is
described by Lemma 3. By an explicit computation, we have the next table
which describes the distributions of
 number of flex points.
The second line is the configuration of singularity and the last line
 is the pair of  flex numbers $(\iota_t,\iota_{nt})$.

\vspace{.3cm}
\noindent
\begin{tabular}{c|c|c|c|c|c|c}

1&2&3&4&5&6&7\\
\hline
  \,\vtop{\hbox{$4A_2$} \hbox{$4A_2+A_1$}\hbox{$4A_2+2A_1$}}&\,\vtop{\hbox{
 $A_5+2A_2$}
\hbox{$A_5+2A_2+A_1$}} & \,\vtop{\hbox{$E_6+2A_2$}\hbox{$E_6+2A_2+A_1$}} &\,$E_6+A_5$&\,$2A_5$&
\,\vtop{\hbox{ $A_8+A_2$}\hbox{$A_8+A_2+A_1$}}&\,$A_{11}$\\
\hline
\vtop{\hbox{(1,12)}\hbox{(1,6)}\hbox{(1,0)}}&\vtop{\hbox{(1,10)}
\hbox{(1,4)}}&\vtop{\hbox{(1,6)}\hbox{(1,0)}}&(1,4)&(1,8)&\vtop{\hbox{(1,9)}\hbox{(1,3)}}&(1,8)\\
\end{tabular}
\end{Remark}
{\em Proof.} First recall that the topology of the complement
of the sextics $B_5\cup L_P$ for a flex point $P$ of torus type
and non-torus type are different. They can be distinguished by 
Alexander polynomial (\cite{OkaAtlas}). Therefore
to show the assertion about the positivity  $\iota_{nt}>0$, it is enough to
check the assertion by some quintic $B_5$. Examples will be given in the
next subsection.  Secondly, the irreducibility of  the configuration space $\cM(\Xi;5)$  of quintics $f_5(x,y)=0$ with  singularities
$\Xi=[4A_2,\,2A_1],\, [E_6,\,2A_2,\,A_1]$ are easily proves as follows.
For $\Xi=[4A_2,\,2A_1]$, the dual curves of quintics in 
this configuration space are quartics with configuration $[A_2,2A_1]$. As
the irreducibility of the configuration space 
$\cM([A_2,2A_1];4)$ is easy to be checked, the irreducibility of
$\cM([4A_2,2A_1];5)$  follows. 
Take $\Xi=[E_6,2A_2,A_1]$.
For a quintic $B_5$ with $\Si(B_5)=\Xi$, the dual curve $B_5^*$ is again a
quartic with $[A_2,\,2A_1]$ (thus  mapped into the same configuration space with the dual  of quintics with
$[4A_2,2A_1]$).
Hoever we can not apply the same argument. The reason is that
 the dual curve $B_5^*$ is
not generic in the configuration space $\cM([A_2,2A_1];4)$: the quartic
$B_5^*$
has not 4 flexes but three flexes, one flex of flex order 4 (=dual of $E_6$) and 2 flexes of flex order
3 (i.e., dual of $2A_2$). Thus we need another argument.
 Note that any three singular points can not be
colinear on
$B_5$ by B\'ezout theorem. We can consider the slice condition:

$(\star)$:  $E_6$ is at $(-1,0)$ and two $A_2$ are at $(0,1),\, (0,-1)$ and one
$A_1$ at $(1,0)$. 

\noindent
It is easy to compute that a Zariski open subset of
this slice has the normal form:
\begin{multline*}
h := \mathit{e_1}^{3} - 4\,y^{3}\,x^{2} - 2\,\mathit{e_1}^{3}\,x^{3
} + x^{4}\,\mathit{e_1}^{3} + x\,\mathit{e_1}^{3} + 3\,y^{5}\,
\mathit{e_1}^{2} + \mathit{e_1}^{3}\,y^{4} + 3\,y\,\mathit{e_1}^{2}
 + 10\,\mathit{e_1}^{3}\,y^{2}\,x^{3} \\
\mbox{} + 18\,\mathit{e_1}^{3}\,y^{2}\,x^{2} - 12\,\mathit{e_1}\,y
^{2}\,x^{2} + 6\,\mathit{e_1}^{3}\,y^{2}\,x - 12\,\mathit{e_1}\,y^{
2}\,x - 9\,y\,x^{4}\,\mathit{e_1}^{2} - 12\,y\,\mathit{e_1}^{2}\,x
^{3} \\
\mbox{} + 12\,y\,x\,\mathit{e_1}^{2} + 6\,y\,x^{2}\,\mathit{e_1}^{2
} - 6\,y^{3}\,x^{2}\,\mathit{e_1}^{2} - 6\,y^{3}\,\mathit{e_1}^{2}
 - 2\,\mathit{e_1}^{3}\,y^{2} + \mathit{e_1}^{3}\,x^{5} - 7\,
\mathit{e_1}^{3}\,y^{4}\,x \\
\mbox{} + 12\,\mathit{e_1}\,y^{4}\,x - 12\,y^{3}\,x\,\mathit{e_1}^{
2} - 2\,x^{2}\,\mathit{e_1}^{3} + 8\,y^{3} - 4\,y - 4\,y^{5} - 8\,
y\,x + 8\,y\,x^{3} + 4\,y\,x^{4} \\
\mbox{} + 8\,y^{3}\,x,\quad e_1\ne 0,\pm  2/\sqrt{3}
\end{multline*}
Thus the irreducibility of $\cM([E_6,2A_2,A_1];5)$ follows. For each of
them we know that the number of flex points is 1 and 
$\cM(\Xi;5)'\ne \emptyset$.
On the other hand, the topology of sextics $B_5\cup L_P$, of torus type and non-torus type,
are distinguished by the Alexander polynomials
$(t^2-t+1)\, (t-1)$ and $(t-1)$ respectively. This implies
$\cM(\Xi;5)'=\cM(\Xi;5)$.
\qed
\subsection{Example for sextics with quintic components}
We gives examples  of sextics with a quintic components.

\noindent
1.
$C=B_5+B_1,\, \Si=[A_5,4A_2, k\,A_1],\,k=2,3$:
First we consider the case $k=2$.  The quintic has $4A_2$ and 
13 flex points.
\begin{multline*}
C:\bigg( ({\displaystyle \frac {28}{153}} \,x + {\displaystyle 
\frac {8}{153}} )\,y^{4} + ({\displaystyle \frac {52}{51}} \,x^{3
} - {\displaystyle \frac {20}{51}} \,x^{2} - {\displaystyle 
\frac {152}{153}} \,x - {\displaystyle \frac {16}{153}} )\,y^{2}
 + x^{5} - {\displaystyle \frac {2}{17}} \,x^{4} - 
{\displaystyle \frac {193}{51}} \,x^{3}\\
\qquad  + {\displaystyle \frac {
116}{51}} \,x^{2} 
\mbox{} + {\displaystyle \frac {124}{153}} \,x + {\displaystyle 
\frac {8}{153}}{\bigg)}
\left({\displaystyle \frac {80}{17}} \,x + 
{\displaystyle \frac {880}{51}}  - {\displaystyle \frac {320}{51}
} \,y\right)
\end{multline*}
The sextics of torus type is obtained by replacing the line component by
$\frac{28}{153}x+\frac{8}{153}$.

Next, we consider the case $k=3$. Th equintic $B_5$ has $4A_2+A_1$.
 The flex which gives a sextic of torus type is $(1,0)$.
\begin{multline*}f_5:=
 - {\displaystyle \frac {385}{16}} \,x^{4}\,y + {\displaystyle 
\frac {3885}{16}} \,x\,y^{2} + {\displaystyle \frac {1897}{128}} 
\,x\,y^{4} + {\displaystyle \frac {345}{16}} \,y\,x^{3} - 
{\displaystyle \frac {441}{4}} \,y^{3}\,x + {\displaystyle 
\frac {529}{4}} \,y\,x^{2} + 73\,y + 72\,y^{3} \\
\mbox{} + {\displaystyle \frac {403}{128}} \,x^{3}\,y^{2} - 
{\displaystyle \frac {16783}{128}} \,x^{2}\,y^{2} - 
{\displaystyle \frac {811}{4}} \,y\,x - {\displaystyle \frac {869
}{64}} \,y^{4} + {\displaystyle \frac {7087}{128}} \,x - 
{\displaystyle \frac {3675}{32}} \,y^{2} + {\displaystyle \frac {
3201}{512}} \,x^{5} \\
\mbox{} + {\displaystyle \frac {601}{32}} \,x^{4} - y^{5} + 
{\displaystyle \frac {313}{8}} \,x^{2}\,y^{3} - {\displaystyle 
\frac {11511}{256}} \,x^{2} - {\displaystyle \frac {997}{64}}  - 
{\displaystyle \frac {10167}{512}} \,x^{3} 
\end{multline*} 
$B_5$ has two obvious flex points: $P:=(1,0)$ and $Q:=( -1520/293, -287/293)$, where $P$ is a flex
 of torus type and $Q$ is a flex of non-torus type. There are 5  other
     flex points whose $x$-coordinates are the solution of
\begin{multline*}
R_1 := 5926214587003\,x^{5} - 32698277751050\,x^{4} + 
69779834665700\,x^{3} \\
\mbox{} - 72918583611000\,x^{2} + 37638730560000\,x - 
7728486400000 =0
\end{multline*}
We can check that the roots of $R_1=0$
corresponds to flexes of non-torus type as follows. (The same argument
applies to other cases.)
Note that any   conics which is passing through 4 $A_2$ of
     $B_5$  are given by
\[
h_2:= y^{2} - {\displaystyle \frac {1}{2}} \,{d_{01}}\,y\,x + 
{d_{01}}\,y - {\displaystyle \frac {1}{2}} \,x^{2} + 
{\displaystyle \frac {1}{19}} \,x^{2}\,{d_{01}} - 
{\displaystyle \frac {5}{2}} \,x - {\displaystyle \frac {10}{19}
} \,{d_{01}}\,x + 3 + {\displaystyle \frac {16}{19}} \,
{d_{01}}
\]
Thus if there is a flex $P(a,b) $ of torus type (so $R_1(a)=0$), there is a cubic form $h_3(x,y)$
such that the sextic $C=B_5\cup L_P$ is described as
 $C: =\{h_3^2+h_2^3=0\}$.
On the other hand, put  $S_2(x,d_{01})$ be the polynomial of degree 2 in
 $x$
defined by  $S_2(x,d_{01})=R(h_2,f_5,y)/P(x)^2$ where $R(h_2,f_5,y)$ is
the resultant of $h_2$ and $f_5$ in $y$ and $P(x)=0$ is the defining polynomial for the x-coordinates of 4 $A_2$.
Then $S_2$ must
be
 $ c\,(x-a)^2$  for some $c\ne 0$. Let $b_1(d_{01})$ be the discriminant polynomial  of $S_{2}$
in $x$ 
and let $b_2(d_{01})$
be the resultant of $S_{2}(x,d_{01})$ and $R_1(x)$
in $x$. Thus we obtain two polynomials
$b_1(d_{01}),\, b_2(d_{01})$
of the parameter $d_{01}$ which must have a common root:
  We can check that  $b_1(d_{01})=b_2(d_{01})=0$
has no common root in $d_{01}$.

\noindent
2. 
 $C=B_5+B_1$, $\Si=[2A_5,2A_2,j\,A_1],\,j=2,3$. The quintic $B_5$ has $A_5+2A_2$.
 \begin{multline*}
j=2,\,[2A_5,2A_2,2A_1]:\,
 \bigg(  - {\displaystyle \frac {16145}{1024}} \,y^{3} - 
 {\displaystyle \frac {93}{64}} \,y^{5} - {\displaystyle \frac {
 877727}{8192}} \,y^{2} - {\displaystyle \frac {110055}{1024}} \,y
  - {\displaystyle \frac {329525}{16384}} \,x \\+ {\displaystyle 
 \frac {11025}{16384}} \,x^{3} - {\displaystyle \frac {543975}{
 16384}} \,x^{2} 
 \mbox{} + {\displaystyle \frac {1751733}{16384}} \,y^{2}\,x^{2}
  + {\displaystyle \frac {79625}{4096}} \,x^{5} - {\displaystyle 
 \frac {235529}{16384}} \,y^{4} + {\displaystyle \frac {1999101}{
 16384}} \,y^{2}\,x^{3}\\ - {\displaystyle \frac {100809}{8192}} \,y
 ^{3}\,x 
 \mbox{} + {\displaystyle \frac {1199495}{8192}} \,y\,x^{3} + 
 {\displaystyle \frac {289995}{4096}} \,y\,x^{2} - {\displaystyle 
 \frac {1199495}{8192}} \,y\,x + {\displaystyle \frac {150225}{
 4096}} \,y\,x^{4}\\ - {\displaystyle \frac {498845}{4096}} \,y^{2}
 \,x 
 \mbox{} - {\displaystyle \frac {275703}{16384}} \,y^{4}\,x + 
 {\displaystyle \frac {23889}{8192}} \,y^{3}\,x^{2} + 
 {\displaystyle \frac {18625}{512}} \,x^{4} - {\displaystyle 
 \frac {52025}{16384}} \bigg)(y+1)
 \end{multline*}
A sextic of torus type is give by replacing the line component by
 $x-1=0$.

 \begin{multline*}
j=3,\, [2A_5,2A_2,3A_1]:\, \bigg({\displaystyle \frac {2}{7}}  + x^{5} - {\displaystyle 
 \frac {4}{7}} \,x^{2} - 2\,x^{3} + {\displaystyle \frac {2}{7}} 
 \,x^{4} + x - {\displaystyle \frac {2}{7}} \,y^{2}\,x^{3} + 
 {\displaystyle \frac {12}{7}} \,x^{2}\,y^{2} - {\displaystyle 
 \frac {1}{7}} \,x\,y^{4} \\
- {\displaystyle \frac {6}{7}} \,y^{2}\,
 x - {\displaystyle \frac {4}{7}} \,y^{2} + {\displaystyle \frac {
 2}{7}} \,y^{4}\bigg)\bigg(
 {\displaystyle \frac {44064}{34157767}} \,y\,\sqrt{ - 963 + 1182
 \,\sqrt{6}} 
+ {\displaystyle \frac {16521840}{34157767}}  - 
 {\displaystyle \frac {7198560}{34157767}} \,\sqrt{6} 
 \mbox{} \\- {\displaystyle \frac {28320}{34157767}} \,y\,\sqrt{ - 
 963 + 1182\,\sqrt{6}}\,\sqrt{6} + {\displaystyle \frac {7328592}{
 34157767}} \,\sqrt{6}\,x - {\displaystyle \frac {2704104}{4879681
 }} \,x{\bigg)}
 \end{multline*}
The quintic $B_5$ has $A_5+2A_2+A_1$
and  5 flex points and among them, there exists  a unique flex of torus type.
The tangent line  at this flex of torus type is given  by $2-x=0$.

\noindent
3. A sextic $C=B_5+B_1$ with $\Si(C)=[E_6,A_5,2A_2,2A_1]$ is given by 
\begin{multline*}
f:=
\big({\displaystyle \frac {53}{141}} \,x + {\displaystyle \frac {3}{
47}} \,y + y^{5} - {\displaystyle \frac {50}{47}} \,y^{3} + 
{\displaystyle \frac {4}{47}} \,y^{2}\,x - {\displaystyle \frac {
769}{141}} \,y^{4}\,x - {\displaystyle \frac {614}{141}} \,y^{2}
\,x^{3} + 2\,y^{2} + {\displaystyle \frac {53}{141}} \,x^{5} + 
{\displaystyle \frac {56}{141}} \,y\,x^{2} \\
\mbox{} - {\displaystyle \frac {10}{3}} \,y\,x + {\displaystyle 
\frac {1174}{141}} \,y^{3}\,x + {\displaystyle \frac {1256}{141}
} \,y^{3}\,x^{2} + {\displaystyle \frac {10}{3}} \,x^{3}\,y - 
{\displaystyle \frac {69}{47}} \,y^{4} - {\displaystyle \frac {25
}{47}} \,x^{4} + {\displaystyle \frac {50}{47}} \,x^{2} - 
{\displaystyle \frac {106}{141}} \,x^{3} \\
\mbox{} - {\displaystyle \frac {1462}{141}} \,y^{2}\,x^{2} - 
{\displaystyle \frac {65}{141}} \,y\,x^{4} - {\displaystyle 
\frac {25}{47}} \big)\big(y + 1 - {\displaystyle \frac {8}{3}} \,x\big) 
\end{multline*}
The quintic has 7 flex points and there is a unique one among them
which is  of torus type at 
$( \frac {2400}{1357}, \frac {357}{1357})$.

\noindent
4. $C=B_5+B_1$ with $[E_6,2A_5,2A_1]$. The quintic $B_5$ has $E_6+A_5$
and it has 5 flexes. Among them, there is a unique flex 
 of torus type. A sextic of non-torus type:
\begin{multline*}
{f} :=( 4451 + 9742\,y^{2}\,x + 4639\,y^{4}\,x - 9501\,y
 - 14381\,x - 423\,y^{5}\,\sqrt{33} - 351\,\sqrt{33} \\
\mbox{} - 16343\,x^{3} + 6546\,y^{3} - 8005\,y^{4} + 3554\,y^{2}
 + 19373\,x^{2} - 19373\,y^{2}\,x^{2} + 9836\,x^{4} \\
\mbox{} + 2955\,y^{5} + 10266\,y\,x^{3} - 2936\,x^{5} - 19020\,y
^{3}\,x + 19020\,y\,x - 14661\,y\,x^{2} \\
\mbox{} + 14661\,y^{3}\,x^{2} + 1521\,\sqrt{33}\,y - 1098\,y^{3}
\,\sqrt{33} + 1593\,y^{4}\,\sqrt{33} + 756\,x^{4}\,\sqrt{33} \\
\mbox{} + 1215\,x^{2}\,\sqrt{33} - 1242\,y^{2}\,\sqrt{33} - 1917
\,x^{3}\,\sqrt{33} + 297\,x\,\sqrt{33} + 999\,y^{2}\,x^{3}\,
\sqrt{33} \\
\mbox{} + 36\,y\,x^{4}\,\sqrt{33} + 1458\,y^{2}\,x\,\sqrt{33} + 
918\,y\,x^{3}\,\sqrt{33} + 2052\,y^{3}\,x\,\sqrt{33} - 423\,y\,x
^{2}\,\sqrt{33} \\
\mbox{} - 2052\,y\,x\,\sqrt{33} + 423\,y^{3}\,x^{2}\,\sqrt{33} - 
1215\,y^{2}\,x^{2}\,\sqrt{33} - 1755\,y^{4}\,x\,\sqrt{33} + 6077
\,y^{2}\,x^{3} \\
\mbox{} - 5124\,y\,x^{4} )(y+1)
\end{multline*}
A sextic torus type is given by replacing $y+1$ by the flex tangent at 
$(\alpha,\beta)$ where
\[
 \alpha  := {\displaystyle \frac {1476423}{6805087}}  + 
{\displaystyle \frac {176748}{6805087}} \,\sqrt{33},\,
\beta  := {\displaystyle \frac {1469468}{6805087}}  - 
{\displaystyle \frac {931392}{6805087}} \,\sqrt{33}
\]

\noindent
5.  $C=B_5+B_1$ with $[3A_5,2A_1]$.
 \begin{multline*}
  - {\displaystyle \frac {12}{6279955}} ( - 2516 + 27\,\sqrt{69})(
  - 2516 - 10064\,x^{4} - 45\,x\,y^{4}\,\sqrt{69} + 90\,x\,y^{2}\,
 \sqrt{69} \\
 \mbox{} + 108\,\sqrt{69}\,y^{3}\,x^{2} - 108\,\sqrt{69}\,y\,x^{2}
  + 180\,x^{5}\,\sqrt{69} + 54\,y^{2}\,\sqrt{69} - 45\,x\,\sqrt{69
 } - 27\,y^{4}\,\sqrt{69} \\
 \mbox{} - 108\,x^{4}\,\sqrt{69} - 27\,\sqrt{69} + 10060\,y\,x^{4}
  - 5030\,y^{3} - 10064\,y\,x^{2} - 2516\,y^{4} \\
 \mbox{} + 10064\,y^{3}\,x^{2} + 2515\,y^{5} + 10060\,x^{2} - 840
 \,x\,y^{4} + 1680\,x\,y^{2} + 5032\,y^{2} - 10060\,y^{2}\,x^{2}
  \\
 \mbox{} + 3360\,x^{5} + 2515\,y - 840\,x) \\
 \left( ({\displaystyle \frac {31104}{2497}}  - {\displaystyle \frac {
 1620}{2497}} \,\sqrt{69})\,y + {\displaystyle \frac {108}{2497}} 
 \,(52\,\sqrt{69} - 499)\,(x - 1)\right) 
 \end{multline*}
 The flex point which gives a sextic of torus type is 
 \[(\alpha,\beta),\,
 \alpha  := {\displaystyle \frac {957138004}{22902646825}}  + 
 {\displaystyle \frac {2339358408}{22902646825}} \,\sqrt{69},\,
 \beta  := {\displaystyle \frac {540908244}{916105873}}  - 
 {\displaystyle \frac {52210443}{916105873}} \,\sqrt{69}
 \]

\vspace{.2cm}\noindent
6.  a. A sextic of torus type with $[A_8,A_5,A_2,2A_1]$ with line component
    is given by:
\begin{multline*}
f:=\left (- 60\,y^{2} + 60\,y - x^{2}\right)^{3} \\
\mbox{} + \bigg( - {\displaystyle \frac {81}{25}} \,y^{3} + ( - 
{\displaystyle \frac {6849}{25}} \,x - {\displaystyle \frac {9639
}{25}} )\,y^{2} + ({\displaystyle \frac {162}{25}} \,x^{2} + 
{\displaystyle \frac {6849}{25}} \,x + {\displaystyle \frac {1944
}{5}} )\,y + x^{3} - {\displaystyle \frac {162}{25}} \,x^{2}\bigg)^{2}
\end{multline*}
and the line component is defined by $y-1=0$.
It has 10 flexes and the flex at $(0,1)$ gives a flex of torus type.
Other  flex tangent lines give a sextic of non-torus type. For example
\begin{multline*}
 {f_6} := ({\displaystyle \frac {324}{25}} \,x^{5} - 
{\displaystyle \frac {26244}{625}} \,x^{4} - {\displaystyle 
\frac {2223126}{625}} \,y^{2}\,x^{3} + {\displaystyle \frac {
40132557}{625}} \,y^{3}\,x^{2} - {\displaystyle \frac {428706}{
625}} \,y\,x^{4} \\
\mbox{} - {\displaystyle \frac {43281837}{625}} \,y^{2}\,x^{2} - 
{\displaystyle \frac {134993439}{625}} \,y^{5} + {\displaystyle 
\frac {271568079}{625}} \,y^{4} - {\displaystyle \frac {8419248}{
125}} \,y^{3} - {\displaystyle \frac {3779136}{25}} \,y^{2} \\
\mbox{} + {\displaystyle \frac {629856}{125}} \,y\,x^{2} + 
{\displaystyle \frac {1733076}{625}} \,y\,x^{3} - {\displaystyle 
\frac {26628912}{125}} \,y^{2}\,x + {\displaystyle \frac {
132035022}{625}} \,y^{3}\,x + {\displaystyle \frac {1109538}{625}
} \,y^{4}\,x) \\
( - {\displaystyle \frac {64039734}{25}} \,y + {\displaystyle 
\frac {18974736}{5}}  + {\displaystyle \frac {33205788}{25}} \,x)
\end{multline*}

\noindent
b. A sextic of torus type with $[A_8,A_5,A_2,3A_1]$ and with component type
$B_5+B_1$ is given by 
\begin{multline*}
f:=
({\displaystyle \frac {15}{2}} \,y^{2} - {\displaystyle \frac {15
}{2}} \,y - 16\,x^{2})^{3} \\
\mbox{} + ( - {\displaystyle \frac {455}{24}} \,y^{3} + ( - 
{\displaystyle \frac {80}{3}} \,x + {\displaystyle \frac {245}{6}
} )\,y^{2} + ({\displaystyle \frac {140}{3}} \,x^{2} + 
{\displaystyle \frac {80}{3}} \,x - {\displaystyle \frac {175}{8}
} )\,y + 64\,x^{3} - {\displaystyle \frac {140}{3}} \,x^{2})^{2}
\end{multline*}
It has the line component $y-1=0$ and the quintic has 4 flex points
among which only the flex (0,1) is of torus type. An example of sextic
of non-torus type is given by
\begin{multline*}
{f_6} := ({\displaystyle \frac {17920}{3}} \,x^{5} - 
{\displaystyle \frac {19600}{9}} \,x^{4} + {\displaystyle \frac {
3500}{3}} \,y^{2}\,x - {\displaystyle \frac {6125}{3}} \,y\,x^{2}
 + {\displaystyle \frac {47600}{9}} \,x^{3}\,y - {\displaystyle 
\frac {30625}{64}} \,y^{2} + {\displaystyle \frac {332125}{192}} 
\,y^{3} \\
\mbox{} - {\displaystyle \frac {11275}{3}} \,y^{3}\,x^{2} + 5800
\,y^{2}\,x^{2} - {\displaystyle \frac {19600}{9}} \,x\,y^{3} + 
{\displaystyle \frac {9100}{9}} \,x\,y^{4} - {\displaystyle 
\frac {44240}{9}} \,x^{3}\,y^{2} + {\displaystyle \frac {40720}{9
}} \,y\,x^{4} \\
\mbox{} - {\displaystyle \frac {1170775}{576}} \,y^{4} + 
{\displaystyle \frac {450025}{576}} \,y^{5})( - {\displaystyle 
\frac {2222000000}{255584169}} \,y + {\displaystyle \frac {
2156000000}{255584169}}  + {\displaystyle \frac {492800000}{
85194723}} \,x) 
\end{multline*}

\vspace{.2cm}
\noindent
7. A quintic with $A_{11}$ has 9 flex points, among which
there exists a unique flex of torus type.
In the following example, our quintic has a flex of torus type at
$(0,1)$
so the the sextic of torus type is given by
\begin{multline*}
f:=\left( - {\displaystyle \frac {28}{25}} \,y^{2} + {\displaystyle 
\frac {28}{25}} \,y - x^{2}\right)^{3}\\
 +\left ({\displaystyle \frac {511}{
100}} \,y^{3} + ( - {\displaystyle \frac {28}{25}} \,x - 
{\displaystyle \frac {7}{500}} )\,y^{2} + ( - {\displaystyle 
\frac {91}{20}} \,x^{2} + {\displaystyle \frac {28}{25}} \,x - 
{\displaystyle \frac {637}{125}} )\,y - x^{3} + {\displaystyle 
\frac {91}{20}} \,x^{2}\right)^{2}
\end{multline*}
and the line component is $y-1=0$. $B_5$ has 8 flex of non-torus
type. We can take one at $(1,-1) $ so that a sextic of non-torus type is
given by
\begin{multline*}
{f_6} := \big({\displaystyle \frac {6176793}{250000}} \,y^{5}
 - {\displaystyle \frac {7154}{625}} \,y^{4}\,x + {\displaystyle 
\frac {7194719}{250000}} \,y^{4} - {\displaystyle \frac {245049}{
5000}} \,y^{3}\,x^{2} - {\displaystyle \frac {859901}{31250}} \,y
^{3} + {\displaystyle \frac {98}{3125}} \,y^{3}\,x \\
\mbox{} + {\displaystyle \frac {35672}{3125}} \,y^{2}\,x - 
{\displaystyle \frac {405769}{15625}} \,y^{2} + {\displaystyle 
\frac {13181}{5000}} \,y^{2}\,x^{2} - {\displaystyle \frac {7}{
250}} \,y^{2}\,x^{3} - {\displaystyle \frac {2548}{125}} \,y\,x^{
3} + {\displaystyle \frac {57967}{1250}} \,y\,x^{2} \\
\mbox{} + {\displaystyle \frac {7833}{400}} \,y\,x^{4} - 
{\displaystyle \frac {8281}{400}} \,x^{4} + {\displaystyle 
\frac {91}{10}} \,x^{5})({\displaystyle \frac {3306744}{15625}} 
\,y + {\displaystyle \frac {1928934}{15625}}  + {\displaystyle 
\frac {275562}{3125}} \,x) 
\end{multline*}

\section{Configuration coming from quartic geometry}
\subsection{Configuration coming from quartic flex geometry}
We consider the sextics  $C$ with component type $B_4+B_1+B_1'$.
The corresponding possible configurations are

(a) 
$\Si(C)=[3A_5+3A_1]$ and $\Si(B_4)=[A_5]$ or
(b)  $ \Si(C)=[2A_5+2A_2+kA_1],\, k=3,4$ and  $\Si(B_4)=[2A_2+(k-3)A_1]$ or
(c) $\Si(C)=[E_6+2A_5+3A_1]$ and $\Si(B_4)=[E_6]$.

\vspace{.3cm}
 Let $P,Q$ be two flex points on $B_4$ and let $L_P,\, L_Q$ be the
 flex tangents.
We say that a pair of flex points $\{P,Q\}$ are a {\em flex pair of torus
type}
if the sextic $B_4\cup L_P\cup L_Q$ is a sextic of torus type.


\begin{Theorem} {\bf Case (a)  $\Si(C)=[3A_5+3A_1]$}. The quartic $B_4$ has one $A_5$ and 
6 flex points and two line components are flex
tangent lines.
There  exist two flex pairs of torus type.
 The other choices give
 sextics of non-torus type.

{\bf Case (b) $ \Si(C)=[2A_5+2A_2+kA_1],\, k=3,4$ }. The quartic $B_4$ has $2A_2$ or $2A_2+A_1$
according to  $k=3$ or $4$ and  $B_4$ has 8 or 2 flex points
 respectively.
For the case, $k=3$, there are both flex pairs of torus type
and of non-torus type.
For $k=4$, the choice of $\{P,Q\}$ is unique and it is a pair of torus
 type.

{\bf Case (c) $\Si(C)=[E_6+2A_5+3A_1]$ }. $B_4$ has two flexes  and  they gives a pair of torus type.
Thus there is no sextic of non-torus type with $E_6+2A_5+3A_2$ with
 component type $B_4+B_1+B_1'$.
\end{Theorem}
{\em Proof.}
As the configuration spaces of quartics with one $A_5$, or $2A_2$ or $2A_2+A_1$ 
or $E_6$ 
 are  connected, it is enough to check the assertion by an example.

For the non-existence, note that a quartic with 
  $\Si(B_4)=2A_2+A_1$ or $\Si(B_4)=[E_6]$ has exactly 2 flexes. Thus the existence of
sextic of torus type $B_4+B_1+B_1'$ with $[2A_5+2A_2+4A_1]$ 
or $[E_6+2A_5+3A_1]$ implies that
there does not exist sextic of non-torus type with these two configurations\qed

\noindent
{\bf Example} I. We consider the quartic $B_4:=\{g_4=0\}$ with one $A_5$:
\begin{multline*}
{g_4} := 639 - 1350\,x^{2} + 351\,x^{4} + 468\,x^{3} - 108
\,x + 288\,y + 1608\,I\,y^{2}\,x^{2}\,\sqrt{3} + 1452\,I\,y\,x\,
\sqrt{3} \\
\mbox{} - 676\,I\,y\,x^{3}\,\sqrt{3} - 1032\,I\,y^{2}\,x\,\sqrt{3
} - 1452\,I\,y^{3}\,x\,\sqrt{3} - 288\,y^{3} - 918\,y^{2} + 279\,
y^{4} \\
\mbox{} - 648\,y^{3}\,x + 1350\,y^{2}\,x^{2} + 108\,y^{2}\,x - 
936\,y\,x^{3} + 648\,y\,x - 432\,I\,y^{2}\,\sqrt{3} + 776\,I\,y^{
3}\,\sqrt{3} \\
\mbox{} - 1608\,I\,x^{2}\,\sqrt{3} - 152\,I\,\sqrt{3} - 776\,I\,y
\,\sqrt{3} + 1032\,I\,x\,\sqrt{3} + 584\,I\,\sqrt{3}\,y^{4} + 728
\,I\,x^{3}\,\sqrt{3} 
\end{multline*}
$B_4$ has an $A_5$ singularity at $(1,0)$ and 6 flexes at 
\begin{multline*}
 P_1=(0,-1),\,P_2=(0,1),\,P3=( - {\displaystyle \frac {130}{1069}}  + 
{\displaystyle \frac {370}{1069}} \,I\,\sqrt{3}, \,
{\displaystyle \frac {263}{1069}}  + {\displaystyle \frac {156}{
1069}} \,I\,\sqrt{3})\\
P4 := ({\displaystyle \frac {2190}{13333}}  - 
{\displaystyle \frac {4790}{39999}} \,I\,\sqrt{3}, \, - 
{\displaystyle \frac {12671}{13333}}  - {\displaystyle \frac {
4492}{39999}} \,I\,\sqrt{3}),\\
P_5=({\displaystyle \frac {2116}{6841}} \,I\,\sqrt{3}
 + {\displaystyle \frac {632}{20160427}} \,\sqrt{90988874 - 
65462588\,I\,\sqrt{3}} + {\displaystyle \frac {4086}{6841}}  
\mbox{} \\
+ {\displaystyle \frac {498}{20160427}} \,I\,\sqrt{3}\,
\sqrt{90988874 - 65462588\,I\,\sqrt{3}},  \\
 - {\displaystyle \frac {9056}{47887}}  - {\displaystyle \frac {
6133}{47887}} \,I\,\sqrt{3} + {\displaystyle \frac {4}{47887}} \,
\sqrt{90988874 - 65462588\,I\,\sqrt{3}})\\
P_6= ( - {\displaystyle \frac {632}{20160427}} \,\sqrt{
90988874 - 65462588\,I\,\sqrt{3}} + {\displaystyle \frac {4086}{
6841}}  \\
\mbox{} - {\displaystyle \frac {498}{20160427}} \,I\,\sqrt{3}\,
\sqrt{90988874 - 65462588\,I\,\sqrt{3}} + {\displaystyle \frac {
2116}{6841}} \,I\,\sqrt{3},  \\
 - {\displaystyle \frac {9056}{47887}}  - {\displaystyle \frac {
6133}{47887}} \,I\,\sqrt{3} - {\displaystyle \frac {4}{47887}} \,
\sqrt{90988874 - 65462588\,I\,\sqrt{3}})
\end{multline*}
It is easy to check that $\{P_1,P_3\},\,\{P_2,P_4\},\, \{P_5,P_6\}$ give
sextic of torus type. The other cases give sextics of non-torus type.
For example, a nice sextic of non-torus type is given  by taking the
 tangent lines at $P_1$ and $P_2$:
$B_1+B_1': (y-1)\,(y+1)=0$.

\noindent
II. Now we consider the quartic with $2A_2$ (case (b) with $k=3$).
\begin{multline*}
f_4(x,y):={\displaystyle \frac {254143}{4096}} \,x^{4} - {\displaystyle 
\frac {251}{16}} \,x^{3} + {\displaystyle \frac {11}{32}} \,y\,x
^{3} - {\displaystyle \frac {5893}{2048}} \,y^{2}\,x^{2} - 
{\displaystyle \frac {2761}{1024}} \,y\,x^{2} - {\displaystyle 
\frac {126093}{2048}} \,x^{2} - {\displaystyle \frac {11}{32}} \,
y\,x + {\displaystyle \frac {1}{32}} \,y^{3}\,x \\
\mbox{} + {\displaystyle \frac {251}{16}} \,x + {\displaystyle 
\frac {5893}{2048}} \,y^{2} - {\displaystyle \frac {1957}{4096}} 
 + {\displaystyle \frac {211}{4096}} \,y^{4} + {\displaystyle 
\frac {2761}{1024}} \,y - {\displaystyle \frac {251}{1024}} \,y^{
3}
\end{multline*}
It has 8 flex points and four flexes are explicitly written as
\[
 P_1=(1,0),\,P_2=(-1,0),\,P_3=(0,-1),\, P_4=( \frac {16064}{64025} , \, \frac {
-61977}{64025})
\]
Sextics of torus type are given by taking tangent lines at 
$\{P_1,P_2\},\, \{P_3,P_4\}$. As a sextic of non-torus type, we can take 
the tangent lines at $P_1,\, P_3$ so that the sextics is given by
 adding two lines
$(x - 1)\,( - 4\,y - 4 + 16\,x) =0$.
The
configuration space of quartic with $2A_2+A_1$ is connected. Each quartic
 has two flex points and with two
 tangent lines $B_1,\,B_1'$, $B_4\cup B_1\cup B_1'$ gives a sextic of
 torus type with $[2A_5,2A_2,4A_1]$.
Thus there is no sextic of non-torus type $C=B_4+B_1+B_1'$ with
 configuration
$[2A_5,2A_2,4A_1]$.
\subsection{Conical geometry of quartic} Now we consider the 
configuration with component type $B_4+B_2$ will be considered here.
The corresponding configurations are
\begin{enumerate}
\item  $\Si(C)=[2A_5,2A_2,2A_1],\,[2A_5,2A_2,3A_1]$.
\item $\Si(C)=[3A_5,2A_1]$.
\item $\Si(C)=[E_6,2A_5,2A_1]$.
\item $\Si(C)=[A_{11},2A_2,2A_1],\,[A_{11},2A_2,3A_1]$.
\item $\Si(C)=[A_{11},A_5,2A_1]$.
\end{enumerate}
Zariski pairs with the above configurations with fixed component type
$B_4+B_2$
can not be explained by the flex geometry. 

We have to generalize the notion of flex points.
Let $B$ be a given irreducible plane curve of degree
 $d$.
Let $\Phi$ be a linear system of conics and let $\alpha$ be the
 dimension of $\Phi$. For a general smooth point
$P\in B$, the maximal intersection number of $I(B,B_2;P)$ for $B_2\in\Phi$
is $\alpha$. We say $P$ is a {\em  conical flex point with respect to $\Phi$}
if the intersection number $I(B,B_2;P)\ge \alpha+1$.
If $\dim\,\Phi=5$ (so $\Phi$ is the family of all conics), we say
simply that
$P$ is  a {\em conical flex point.}

\vspace{.5cm}\noindent
1. Let us consider the case  $\Si(C)=[2A_5,2A_2,2A_1],\,[2A_5,2A_2,3A_1]$. We consider 
first a sextic of torus type $C=\{f_2^3+f_3^2=0\}$
which decomposes into a quartic $B_4$ and  a quartic $B_2$:

\begin{multline*}
f(x,y):=\left(y^{2} - 2 + 2\,x^{2}\right)^{3} +\\
\left ( - 25\,y^{3} + (13\,x - 23)\,y^{2}
 + ( - 26\,x^{2} + 26)\,y + 13\,x^{3} - 23\,x^{2} - 13\,x + 23\right)^{
2}
\\
=(y^{2} + x^{2} - 1)(177\,x^{4} - 598\,x^{3} - 676\,y\,x^{3} + 344
\,x^{2} + 849\,x^{2}\,y^{2} + 1196\,y\,x^{2} - 650\,y^{3}\,x \\
\mbox{} + 598\,x + 676\,y\,x - 598\,x\,y^{2} - 521 - 151\,y^{2}
 + 1150\,y^{3} - 1196\,y + 626\,y^{4}) 
\end{multline*}
Our quartic is defined by
\begin{multline*}
g_4(x,y)=(177\,x^{4} - 598\,x^{3} - 676\,y\,x^{3} + 344
\,x^{2} + 849\,x^{2}\,y^{2} + 1196\,y\,x^{2} - 650\,y^{3}\,x \\
\mbox{} + 598\,x + 676\,y\,x - 598\,x\,y^{2} - 521 - 151\,y^{2}
 + 1150\,y^{3} - 1196\,y + 626\,y^{4}) =0
\end{multline*}
Note that the singularities $B_4$ are  $2\, A_2$.
The intersection $B_2\cap B_4$ makes  two $A_5$ at $P:=(1,0)$ and $Q:=(-1,0)$.
We consider the linear system $\Phi$ of conics of dimension 2  which
are defined by the
 conics $C_2:=\{h_2(x,y)=0\}$ such that 
$I(C_2,B_4;P)\ge 3$.
Then  we consider the
conical
flex points $R=(a,b)\in B_4$ with respect to $\Phi$, which is described
by the
condition
$\exists h_2\in \Phi$ such that $I(h_2,B_4;R)\ge 3$. We found that there 
are $11$ conical flex points. Two of them can be explicitly given as
$S_1=Q$ and $S_2:=(0,-1)$. The corresponding conics are given as
\begin{eqnarray*}
&h_2(x,y)=x^2+y^2-1,\, \{h_2=0\}\cap B_4\supset \{P,Q\}\\
&g_2(x,y)=(y^{2} + 2\,y\,x - 2\,y - x^{2} + 4\,x - 3),\, \{g_2=0\}\cap B_4\supset \{P,S_2\}
\end{eqnarray*}
We can easily check that the sextic
\begin{multline*}
(y^{2} + 2\,y\,x - 2\,y - x^{2} + 4\,x - 3)(177\,x^{4} - 598\,x^{
3} - 676\,y\,x^{3} + 344\,x^{2} + 849\,x^{2}\,y^{2} \\
\mbox{} + 1196\,y\,x^{2} - 650\,x\,y^{3} + 598\,x + 676\,y\,x - 
598\,x\,y^{2} - 521 - 151\,y^{2} + 1150\,y^{3} \\
\mbox{} - 1196\,y + 626\,y^{4}) =0
\end{multline*}
is not of torus type. Thus 
$B_4\cup\{h_2(x,y)=0\}$ and $B_4\cup \{g_2(x,y)=0\}$
is a Zariski pair.

\vspace{.3cm}\noindent
 Similarly the case $[2A_5,2A_2,3A_1]$ can be treated in the same way.
We start from a sextic of torus type:
\begin{multline*}
f(x,y)=\left( - {\displaystyle \frac {49}{64}} \,y^{2} - {\displaystyle 
\frac {15}{64}}  + {\displaystyle \frac {15}{64}} \,x^{2}\right)^{3} + 
 \\
\left( - {\displaystyle \frac {131}{256}} \,y^{3} + ({\displaystyle 
\frac {729}{512}} \,x - {\displaystyle \frac {297}{256}} )\,y^{2}
 + ( - {\displaystyle \frac {387}{256}} \,x^{2} + {\displaystyle 
\frac {387}{256}} )\,y + {\displaystyle \frac {729}{512}} \,x^{3}
 - {\displaystyle \frac {297}{256}} \,x^{2} - {\displaystyle 
\frac {729}{512}} \,x + {\displaystyle \frac {297}{256}} \right)^{2}\\
={\frac {27}{262144}}\, \left( {y}^{2}+{x}^{2}-1 \right)  ( 19808\,{x}^{4}-41796\,y{x}^{3}-32076\,{x}^{3}+40521\,{y}^{2}{x}^{2}-6865\,{x}^{2}+34056\,{x}^{2}y\\
-32076\,x{y}^{2}\mbox{}+32076\,x+41796\,xy-14148\,x{y}^{3}-7770\,{y}^{2}-34056\,y\\
-1815\,{y}^{4}+11528\,{y}^{3}-12943
\mbox{} ) 
\end{multline*}
The quartic $B_4$ has two $2A_2+A_1$ and we consider the linear system
$\Phi$ of conics intersecting $B_4$ at $P=(1,0)$ with intersection number 3.
We find that there exist 5 conical flex points with respect to $\Phi$,
and among them we have two explicit ones: $Q=(-1,0)$ and $(0,-1)$.
We see that the conic corresponding to $(0,-1)$ gives a Zariski partner 
sextic $f_6=0$ to $C=\{f=0\}$.
\begin{multline*}
\mathit{f_6} := \left( 5\,{y}^{2}-64\,xy+64\,y+69\,{x}^{2}-128\,x+59 \right) 
\mbox{} ( 19808\,{x}^{4}-41796\,y{x}^{3}-32076\,{x}^{3}+40521\,{y}^{2}{x}^{2}\\
-6865\,{x}^{2}+34056\,{x}^{2}y-32076\,x{y}^{2}
\mbox{}+32076\,x+41796\,xy-14148\,x{y}^{3}-7770\,{y}^{2}-34056\,y\\
-1815\,{y}^{4}+11528\,{y}^{3}-12943
\mbox{} )
\end{multline*}

\begin{Remark}The calculation of conical flex points are usually very
 heavy. We used maple 7 to compute in the following recipe.
a. First compute the normal form of $h_2\in \Phi$. It contains two
 parameters.
b. Assume $(u,v)\in B_4$. Put $gg_4(x,y):=(x+u,y+v)$ and
 $hh_2(x,y):=h_2(x+u,y+v)$. Consider the maximal contact 
coordinate at $(u,v)$: $\Phi(x)=a_1\,x+a_2\, x^2+a_3\,x^3$ and put
$GG4(x):=gg_4(x,\Phi(x))$ and $HH_2(x):=hh_2(x,\Phi(x))$.
 Our assumption
 implies that $\Coeff(GG4,x,1)=\Coeff(GG2,x,2)=0$ and
 $\Coeff(HH_2,x,0)=\Coeff(HH_2,x,1)=\Coeff(HH_2,x,2)=0$.
Solve the equations  $\Coeff(GG4,x,1)=\Coeff(GG2,x,2)=0$ in $a_1,\,
 a_2$.
Then solve the equations
$\Coeff(HH_2,x,0)=\Coeff(HH_2,x,1)=0$ in the remaining parameters of the
 linear system. Then we get two equations in $u,v$:
\[
 g_4(u,v)=\Coeff(HH_2,x,2)=0
\]
c. Use the resultant computation to solve the above equations to obtain
 the possibility of conical flex points.
\end{Remark}

\noindent
2. $\Sigma(C)=[3A_5,2A_1]$ with $B_4+B_2$:
\begin{multline*}
 B_4:\, g_4(x,y)=(6\,y\,x + {\displaystyle \frac {1710}{91}} \,y\,
x^{2} - {\displaystyle \frac {1466}{91}} \,x^{3}\,y - 
{\displaystyle \frac {1992}{91}} \,y^{2}\,x - 6\,y^{3}\,x + 
{\displaystyle \frac {790}{91}} \,y^{3} + y^{4} - {\displaystyle 
\frac {4904}{91}} \,x^{3} + {\displaystyle \frac {1992}{91}} \,x
 \\
\mbox{} - {\displaystyle \frac {790}{91}} \,y + {\displaystyle 
\frac {939}{91}} \,y^{2}\,x^{2} + {\displaystyle \frac {1161}{91}
} \,y^{2} + {\displaystyle \frac {1968}{91}} \,x^{2} - 
{\displaystyle \frac {1252}{91}}  + {\displaystyle \frac {2196}{
91}} \,x^{4})
\end{multline*}
It has an $A_5$ at $P:=(1,0)$. We consider the linear system $\Phi$ of conics
of dimension 2  whose conic are intersecting with $B_4$ at $(0,1)$ with
intersection number 3. We find 14 conical flex points with respect to
$\Phi$ in which two are explicit: $R=(\frac {498727}{500817},\frac
{-266266}{500817})$
and $Q=(0,-1)$. The corresponding conics $f_2=0, \,k_2=0$ intersecting
$B_4$ at $P,Q$ oe $P,R$
are given by the following and they gives sextics of non-torus type and of
torus type respectively.
\begin{eqnarray*}
& f_2(x,y):=y^{2} - 1 + {\displaystyle \frac {171}{79}} \,x^{2}
\\
& k_2(x,y)=(y^{2} + ({\displaystyle \frac {268}{759}} \,x + {\displaystyle 
\frac {4424}{2277}} )\,y + {\displaystyle \frac {85291}{19987}} 
\,x^{2} - {\displaystyle \frac {268}{759}} \,x - {\displaystyle 
\frac {6701}{2277}} )
\end{eqnarray*}
\vspace{.3cm}
\noindent
3. $\Si(C)=[E_6,2A_5,2A_1]$: 
We start the next quartic
\begin{multline*}
B_4:\,y^{4} + ( - {\displaystyle \frac {195}{64}} \,x
 + {\displaystyle \frac {169}{64}} )\,y^{3} + ({\displaystyle 
\frac {105}{32}} \,x^{2} - {\displaystyle \frac {33}{8}} \,x + 
{\displaystyle \frac {27}{32}} )\,y^{2} + ( - {\displaystyle 
\frac {143}{64}} \,x^{3} + {\displaystyle \frac {117}{64}} \,x^{2
} + {\displaystyle \frac {195}{64}} \,x - {\displaystyle \frac {
169}{64}} )\,y \\
\mbox{} + {\displaystyle \frac {45}{32}} \,x^{4} - 
{\displaystyle \frac {19}{8}} \,x^{3} - {\displaystyle \frac {21
}{16}} \,x^{2} + {\displaystyle \frac {33}{8}} \,x - 
{\displaystyle \frac {59}{32}} =0
\end{multline*}
Sextics of non-torus type and of torus type are given by the conics
$B_2,\, B_2'$:
\begin{eqnarray*}
& B_2:\, y^{2} - 1 + {\displaystyle \frac {9}{13}} \,x^{2}=0\\
 &B_2':\,y^{2} + ( - {\displaystyle \frac {770}{1147}} \,x
 + {\displaystyle \frac {156}{1147}} )\,y + {\displaystyle 
\frac {11025}{14911}} \,x^{2} + {\displaystyle \frac {770}{1147}
} \,x - {\displaystyle \frac {1303}{1147}}=0
\end{eqnarray*}
\vspace{.3cm}
\noindent
4. We consider the configurations $[A_{11},2A_2,2A_1],\,
[A_{11},2A_2,3A_1]$.
First we consider two cuspidal quartics with $[2A_2]$:
\begin{eqnarray*}
\mathit{f_4} := 5805\,x^{4} - 2916\,I\,x^{3}\,\sqrt{2} + 3888\,I
\,x^{3}\,y\,\sqrt{2} - 1269\,x^{2}\,y^{2} - 729\,x^{2} - 3834\,x
^{2}\,y \\
\mbox{} - 3888\,I\,x\,y^{2}\,\sqrt{2} + 108\,I\,x\,\sqrt{2}\,y^{3
} + 2916\,I\,x\,y\,\sqrt{2} + 1323\,y^{3} - 1971\,y^{2} - 81\,y^{
4} \\
\mbox{} + 729\,y =0
\end{eqnarray*}
It can makes sextics of torus type and non-torus type with configuration
$[A_{11},2A_2,2A_1]$ with respective
conics:
\begin{eqnarray*}
&f_2(x,y)=\,y-x^2\\
&h_2(x,y)=\,({\displaystyle \frac {1}{8}} \,I\,\sqrt{2}\,y^{
2} + x\,y - {\displaystyle \frac {3}{4}} \,I\,y\,\sqrt{2} - 
{\displaystyle \frac {5}{4}} \,I\,\sqrt{2}\,x^{2} - 3\,x + 
{\displaystyle \frac {5}{8}} \,I\,\sqrt{2})
\end{eqnarray*}
They correspond to the conical flex points $(0,0)$ and (0,1).
The other conical flexes are very heavy to be computed.

Next we consider the configuration $[A_{11},2A_2,3A_1]$ which is
produced by a quartic $B_4$ with $2A_2+A_1$ and a conic $B_2$ with a single tangent at
a conical flex.
\[ B_4:
{\displaystyle \frac {1}{16}} \,y^{4} + {\displaystyle \frac {3}{
4}} \,x\,y^{3} + {\displaystyle \frac {59}{8}} \,y^{2}\,x^{2} - y
^{2}\,x - {\displaystyle \frac {1}{8}} \,y^{2} + {\displaystyle 
\frac {27}{4}} \,y\,x^{3} - 6\,y\,x^{2} - {\displaystyle \frac {3
}{4}} \,y\,x + {\displaystyle \frac {17}{16}} \,x^{4} - x^{3} - 
{\displaystyle \frac {9}{8}} \,x^{2} + x + {\displaystyle \frac {
1}{16}}=0 
\]
We  find three conical flex points
\[
 P1 := (-1/19, 12/19),\quad
\mathit{P2} := ({\displaystyle \frac {-4}{13}} , \,
{\displaystyle \frac {15}{13}} ),\quad
\mathit{P3} := (-1, \,0)
\]
The corresponding conic which are tangent  at the respective conical
flex point $P_i,\,i=1,2,3$
are given by
\begin{eqnarray*}
g_{21}:=y^{2} + ({\displaystyle \frac {1270}{141}} \,x + {\displaystyle 
\frac {10}{141}} )\,y + {\displaystyle \frac {38711}{423}} \,x^{2
} - {\displaystyle \frac {5486}{423}} \,x - {\displaystyle 
\frac {457}{423}}\\
g_{22}:=y^{2} + ({\displaystyle \frac {7462}{2517}} \,x - {\displaystyle 
\frac {5462}{2517}} )\,y + {\displaystyle \frac {13007}{7551}} \,
x^{2} - {\displaystyle \frac {21902}{7551}} \,x + {\displaystyle 
\frac {8831}{7551}}\\
g_{23}:=- {\displaystyle \frac {3}{22}} \,y^{2} - {\displaystyle \frac {
61}{11}} \,y\,x + y - {\displaystyle \frac {73}{66}} \,x^{2} - 
{\displaystyle \frac {1}{33}} \,x + {\displaystyle \frac {71}{66}
}
\end{eqnarray*}
The corresponding sextic $f_4(x,y)\,g_{2\,latex Conic-all-Ri}(x,y)=0$
is of torus type for $i=1$ and of non-torus type for i=2,3.
The torus decomposition of $f_4\, g_{21}$ is given by 
$z_{21}^3+z_{31}^2=0$ where 
\begin{eqnarray*}
&\mathit{z_{21}} := {\displaystyle \frac {1}{423}} \,423^{(2/3)}\,(3
\,y^{2} + 24\,y\,x + 293\,x^{2} - 34\,x - 3)\\
&\mathit{z_{31}} := {\displaystyle \frac {1}{6768}} \,I 
(5 - 189\,x - 3477\,x^{2} + 20045\,x^{3} - 3\,y - 54\,y\,x \\
&\qquad+ 2361
\,y\,x^{2} - 5\,y^{2} + 247\,y^{2}\,x + 3\,y^{3}) 
\sqrt{6768} \end{eqnarray*}

\vspace{.3cm}
\noindent
5. Lastly, we consider the configuration $[A_{11},A_5,2A_1]$ which is
associated to a quartic $B_4$ with an $A_5$ singularity and a conic
$B_2$ tangent at a conical flex point with intersection number 6.
As a quartic, we take:
\begin{multline*}\mathit{f4} := y^{4} + x\,y^{3} + {\displaystyle \frac {7}{15}} 
\,x^{2}\,y^{2} - 2\,x\,y^{2} - 3\,y^{2} + {\displaystyle \frac {2
}{15}} \,x^{3}\,y - {\displaystyle \frac {2}{15}} \,x^{2}\,y + x
\,y + 2\,y + {\displaystyle \frac {4}{75}} \,x^{4}\\ - 
{\displaystyle \frac {2}{15}} \,x^{3} 
\mbox{} - {\displaystyle \frac {1}{3}} \,x^{2}
\end{multline*}
$B_4$ has apparently 26 conical flex points. (The calculation is very
heavy.)
We take four explicit conical flex points:
\[
 P_1:=(0,0),\, P_2:=(-\frac{540}{493},\frac {250}{493}),\,
P_3:=( \frac {-270}{301}, \frac {58}{301}),\,P_4:=(\frac {-270}{193}  , \frac {-50}{193})
\]
After an easy computation, the respective conics are given as
\begin{eqnarray*}
&\mathit{n_{21}} := y^{2} + ( - {\displaystyle \frac {5}{59}} \,x - 
{\displaystyle \frac {50}{59}} )\,y + {\displaystyle \frac {25}{
177}} \,x^{2}\\
&\mathit{n_{22}} := y^{2} + ( - {\displaystyle \frac {10845}{262699}
} \,x - {\displaystyle \frac {273650}{262699}} )\,y + 
{\displaystyle \frac {135675}{262699}} \,x^{2} + {\displaystyle 
\frac {153000}{262699}} \,x + {\displaystyle \frac {70000}{262699
}} \\
&\mathit{n_{23}} := y^{2} + ( - {\displaystyle \frac {16681}{32607}} 
\,x - {\displaystyle \frac {17318}{10869}} )\,y + {\displaystyle 
\frac {30251}{163035}} \,x^{2} - {\displaystyle \frac {1544}{
32607}} \,x - {\displaystyle \frac {112}{10869}} \\
&\mathit{n_{24}} := y^{2} + ({\displaystyle \frac {5225}{1633}} \,x
 + {\displaystyle \frac {350}{71}} )\,y + {\displaystyle \frac {
4225}{1633}} \,x^{2} + {\displaystyle \frac {13000}{1633}} \,x + 
{\displaystyle \frac {10000}{1633}} 
\end{eqnarray*}
Put $f_{6j}(x,y):= f_4(x,y)\, n_{2j}(x,y)$ and $C^{(j)}=\{f_{6j}=0\}$.
It is also easy to see that $C^{(1)},\, C^{(2)}$ are of non-torus type
and $C^{(3)},\,C^{(4)}$ are of torus type. 

\section{Flex geometry of cubic curves}
\subsection{Configurations coming from cubic flex geometry: a cubic
  component and a line component.}
Let us consider first configurations which occurs in sextics which have
at
least a cubic component $B_3$. We divides into the following cases.
\begin{enumerate}
\item $C=B_3+B_3'$.
\begin{enumerate}
\item
$\Si(C)=[A_{17}],\,[A_{17},A_1],\, [A_{17},2A_1]$.
\item $\Si(C)=[A_{11},A_5],\,[A_{11},A_5,A_1],\,[A_{11},A_5,2A_1],\,$.
\item $\Si(C)=[A_{11},2A_2,3A_1]$
\item $\Si(C)=[3A_5],\,[3A_5+A_1],\,[3A_5,2A_1]$.
\end{enumerate}
\item $C=B_3+B_2+B_1$.
\begin{enumerate}
\item $\Si(C)=[A_{11}+A_5+2A_1],\,[A_{11}+A_5+3A_1]$.
\item $\Si(C)=[3A_5+2A_1],\,[3A_5+3A_1]$ .
\end{enumerate}
\item $C=B_3+B_1+B_1'+B_1''$ with configuration
      $[3A_5+3A_1],\,[3A_5+4A_1]$.
\end{enumerate}

We first consider the cases (2) and (3). In these cases, there are one cubic
component
$B_3$ and at least one line component $B_1$.
Recall that the configurations in $(2)$ and $(3)$ occurs as sextics of
linear torus type.
For a  reducible sextic  $C$ which is classified in either (2) or (3),
the necessary and sufficient condition for $C$ to be of torus type is
there exists a line $L$ containing
inner singularities. In the case of  $\Si(C)=[A_{11},A_5]$, 
$L$ is also tangent to the tangent cone of $A_{11}$.
%
%
%
%
We first recall the following basic geometry for cubic curves.
\begin{Proposition}
1. 
A smooth cubic $C$ has 9 flex points. Among $84$ choices of three flex
 points,
12 colinear triples of flexes.

2. A nodal cubic has 3 flex points, and they are colinear.

3. A cuspidal cubic has one flex point.
\end{Proposition}
For the proof of the assertion 1, see Example below.
\begin{Corollary} \label{non-existence1}The configuration $[3A_5+4A_1]$ with components type
$B_3+B_1+B_1'+B_1''$
does not exist as a sextic of non-torus type.
\end{Corollary}
{\em Proof.} The cubic has a node and three line components are flex
tangent
lines at three flex points. 
We know that 
such configuration exists as a sextics of linear torus type
\cite{Reduced}.
As the  configuration  space of one nodal cubics is connected,
every sextics $B_3+B_1+B_1'+B_1''$ is of torus type.\qed

\vspace{.3cm}
\noindent
{(2)  $C=B_3+B_2+B_1$ with $\Si(C)=[3A_5,kA_1],\, k=2,3$.}
In this case, $B_3$ is either smooth or nodal and two intersection points
 $B_3\cap B_2$
generates  $2\,A_5$. The  third $A_5$ is generated by a flex tangent
line $B_1$. 

\begin{Proposition}\label{flex1}
Assume that a cubic $B_3$ 
and a conic $B_2$  are intersecting at two points 
$P,Q$ with respective intersection number 3, producing $2A_5$-singularities.
Then the line passing through $P,\,Q$ intersects $B_3$ at another point,
 say $R$, and  $R$ is a flex point of $B_3$.
\end{Proposition}
This Proposition describes sextics of torus type and non-torus type
with configuration $[3A_5,j\,A_1],\,j=0,1,2$ ((2-b)).

{\em Proof of Proposition \ref{flex1}.} Assume that $P=(0,1)$ and $Q=(0,-1)$ with 
the tanget lines $y=\pm 1$ respectively.  Then by an easy computation,
the cubic $B_3$ is defined by a polynomial
\[
 f_3:=y^{3} + y^{2}\,{a_{12}}\,x - y^{2}\,{a_{00}} + y\,
{a_{21}}\,x^{2} - y + {a_{30}}\,x^{3} - {a_{00}}\,
{a_{21}}\,x^{2} - {a_{12}}\,x + {a_{00}}=0
\]
and the conic is given as 
$y^{2} + {a_{21}}\,x^{2} - 1=0$.
Then $R$ is given as $(0,a_{00})$ and we can easily see that $R$ is a
flex of $B_3$.\qed

By the same calculation, we see that 
\begin{Proposition}\label{non-flex}
Assume that  a cubic $B_3$  and a
 conic $B_2$  are intersecting at one points 
$P$ with intersection multiplicity $6$ producing an $A_{11}$-singularity
Then the tangent line passing at  $P$ intersects another point $R\in
 B_3$ and  $R$ is a flex point of $B_3$.
\end{Proposition}
Proof is similar. Putting $P=(0,0)$ and assuming $y=0$ as the tangent
line,
the cubic is written as 
\begin{multline*}
{f_3} := (y^{3}\,{t_2}^{4} - y^{2}\,x\,{t_3}\,
{a_{01}}\,{t_4}\,{t_2} + y^{2}\,x\,{a_{21}}\,
{t_3}\,{t_2}^{2} + 2\,y^{2}\,x\,{a_{11}}\,
{t_3}^{2}\,{t_2} + 2\,y^{2}\,x\,{a_{01}}\,
{t_3}^{3} \\
\mbox{} - y^{2}\,x\,{a_{11}}\,{t_4}\,{t_2}^{2} - 
y^{2}\,{t_2}^{2}\,{a_{01}}\,{t_4} - y^{2}\,
{t_2}^{3}\,{a_{21}} - y^{2}\,{t_2}^{2}\,{
a_{11}}\,{t_3} + y\,{t_2}^{4}\,{a_{21}}\,x^{2} \\
\mbox{} + y\,{t_2}^{4}\,{a_{11}}\,x + y\,{t_2}^{4
}\,{a_{01}} - x^{3}\,{t_2}^{4}\,{a_{01}}\,{
t_3} - x^{3}\,{t_2}^{5}\,{a_{11}} - {a_{01}}\,
{t_2}^{5}\,x^{2})/{t_2}^{4} 
\end{multline*}
and $R=(0,
x\,{t_2}\,{a_{11}} + {a_{01}}\,{t_2} + x\,
{a_{01}}\,{t_3})$.

\begin{Lemma}
Assume that a conic $B_2$ is tangent to an irreducible curve
 $C$ of degree $d\ge 3$ at a smooth point $P\in C$ so that
 $I(B_2,C;P)\ge 3$. Then $P$ 
is not a flex point of $C$.
\end{Lemma}
{\em Proof.} Let $h_2(x,y)=0$ be a conic equation which defines $B_2$.
In fact, if $P=(a,b)$ is a flex point of $C$, $I(C,B_2;P)\ge 3$ implies that
$C$
is locally parametrized as $y_1(x)=t_1\,x_1+t_3\, x_1^3+\text{(higher
terms)}$ where $(x_1,y_1)=(x-a,y-b)$ assuming the tangent line is not $x-a=0$. 
As $B_2$ does not have any flex, the equation $h_2(x,y)=0$ is solved as 
$y_1=s_1\, x_1+s_2\, x_1^2+\text{(higher terms)}$ with $s_2\ne 0$.
 Thus $I(C,B_2;P)=1$ or $2$ according to $s_1\ne t_1$ or $s_1=t_1$.\qed

First we consider the case (2-b). Then the cubic is either smooth or have
a node. Thus it has at least 3 flexes. As the intersection 
$B_3\cap B_2$
are not flex points, we can find another flex point $S\in B_3$, $S\,\ne\, R$. Taking
flex tangent $L:=T_S \,B_3$ as the line component, the corresponding
sextic is not of torus type. As a cuspidal cubic has a unique flex
points,
we see that a sextic $C=B_3+B_2+B_1$ with $[3A_5+A_2+2A_1]$ does not
appear  as a sextic of non-torus type.

Now we consider (3). Assume that the cubic is smooth. Then there are 9
flex points and $B_1,B_1',B_1''$ are flex  tangents. The sextic is of
torus type
if and only if three flexes are colinear \cite{OkaAtlas}. This case, the
configuration is $[3A_5+3A_1]$.

Finally we consider the case (2-a). In this case, $B_3\cap B_2$ is a
single point $P$ and $I(B_3,B_2;P)=6$ and the intersection singularity
is $A_{11}$. $B_3$ has  at most a node and so it has at least
3 flex points. Taking a line component which is the flex tangent at $S$ other
than $R$, we get a sextic $C=B_3+B_2+B_1$ with $[A_{11}+A_5+2A_1]$ 
or
$[A_{11}+A_5+3A_1]$. \qed

We omit explicit examples  for (2-a) and (2-b) as they can be easily 
obtained from sextic of torus type with the same configuration
(\cite{Reduced}) and
replacing the flex line $B_1$.
We only gives an example of (3).

\vspace{.3cm}
{\bf Example  $B_3+B_1+B_1'+B_1''$, with configuration $[3A_5+3A_1]$}.
In this case, the cubic $B_3$ is smooth and has $9$ flexes and three
lines $B_1,B_1',B_1''$
are the tangent lines at flexes $P_1,\,P_2,\,P_3$ (see below) of $B_3$. 
Three $A_1$ are the intersections of
lines.
 We know that  $B_3\cup L_P\cup L_Q\cup L_R$ is of torus type if and
only if $P,Q,R$ are colinear, where $L_P$ is the tangent line at
the flex point $P$. Let $B_3: f_3(x,y)=0$.
An Example of such a sextic of non-torus type is given by
\begin{eqnarray*}
&f_3(x,y)=y^{3} + ( - ( - 3 + {\displaystyle \frac {1}{2}} \,I\,\sqrt{3})\,
x - 2)\,y^{2} - y + (1 - {\displaystyle \frac {1}{2}} \,I\,\sqrt{
3})\,x^{3} + ( - 3 + {\displaystyle \frac {1}{2}} \,I\,\sqrt{3})
\,x + 2\\
&f(x,y)= (y^{3} + ( - ( - 3 + {\displaystyle \frac {1}{2}} \,I\,
\sqrt{3})\,x - 2)\,y^{2} - y + (1 - {\displaystyle \frac {1}{2}} 
\,I\,\sqrt{3})\,x^{3}
 + ( - 3 + {\displaystyle \frac {1}{2}} \,I\,\sqrt{3})\,x + 2)\\
&(y - 1) 
(y + 1)\,( - I\,\sqrt{3}\,(x - 1) - y) 
\end{eqnarray*}
We can moreover explicitly compute 9 flex points $P_1,\dots, P_9$ as follows.
\begin{multline*}
{P_1} := (1, \,0, \,1),\,
{P_2} := (0, \,1, \,1),\,
{P_3} := (0, \,-1, \,1),\,
{P_4} := (0, \,2, \,1),\,
{P_5} := ({\displaystyle \frac {3}{5}} , \,{\displaystyle 
\frac {4}{5}} , \,1),\,\\
{P_6} := ({\displaystyle \frac {1}{2}}  + {\displaystyle 
\frac {1}{6}} \,I\,\sqrt{3}, \,{\displaystyle \frac {1}{2}}  - 
{\displaystyle \frac {1}{6}} \,I\,\sqrt{3}, \,1),\,
{P_7} := ({\displaystyle \frac {15}{14}}  + {\displaystyle 
\frac {3}{14}} \,I\,\sqrt{3}, \,{\displaystyle \frac {1}{14}}  + 
{\displaystyle \frac {3}{14}} \,I\,\sqrt{3}, \,1),\,\\
{P_8} := ({\displaystyle \frac {21}{38}}  + {\displaystyle 
\frac {9}{38}} \,I\,\sqrt{3}, \,{\displaystyle \frac {31}{38}} 
 - {\displaystyle \frac {3}{38}} \,I\,\sqrt{3}, \,1),\,
{P_9} := ({\displaystyle \frac {33}{62}}  + {\displaystyle 
\frac {3}{62}} \,I\,\sqrt{3}, \,{\displaystyle \frac {37}{62}} 
 + {\displaystyle \frac {9}{62}} \,I\,\sqrt{3}, \,1)
\end{multline*}
Thus by a direct checking, we find the following 12 triples which are
colinear.

\begin{eqnarray*}
 &{\cC_1} := [{P_1}, \,{P_2}, \,{P_6}],\,
 {\cC_2} := [{P_1}, \,{P_3}, \,{P_7}],\,
 {\cC_3} := [{P_1}, \,{P_4}, \,{P_5}]\\
& {\cC_4} := [{P_1}, \,{P_8}, \,{P_9}],\,
 \cC_5 := [{P_2}, \,{P_3}, \,{P_4}],\,
{\cC_6} := [{P_2}, \,{P_5}, \,{P_8}]\\
& {\cC_7} := [{P_3}, \,{P_5}, \,{P_9}],\,
 {\cC_8} := [{P_3}, \,{P_6}, \,{P_8}],\,
 {\cC_9} := [{P_4}, \,{P_6}, \,{P_9}]\\
 &{\cC_{10}} := [{P_4}, \,{P_7}, \,{P_8}],\,
 {\cC_{11} }:= [{P_5}, \,{P_6}, \,{P_7}],\,
 {\cC_{12}} := [{P_2}, \,{P_7}, \,{P_9}]
\end{eqnarray*}
\subsection{Sextics with two cubic components: $C=B_3+B_3'$.}

Now we consider sextics with  two cubic curves $B_3,\,B_3'$.
The possible configurations are
\begin{enumerate}
\item $C=B_3+B_3'$.
\begin{enumerate}
\item
$\Si(C)=[A_{17}],\,[A_{17},A_1],\, [A_{17},2A_1]$.
\item $\Si(C)=[A_{11},A_5],\,[A_{11},A_5,A_1],\,[A_{11},A_5,2A_1],\,$.
\item $\Si(C)=[A_{11},2A_2,3A_1]$
\item $\Si(C)=[3A_5],\,[3A_5+A_1],\,[3A_5,2A_1]$.
\end{enumerate}
\end{enumerate}
First we consider two cubics $B_3,\,B_3'$
 which are
tangent at the origin with intersection number $9$.
Let $f(x,y)=0$ and $f'(x,y)=0$ be the defining polynomials of $B_3$ and
$B_3'$ respectively and we may assume that the tangent line of $B_3$ is
given by $y=0$.
Let $y=\sum_{i=2}^\infty t_i x^i$ be the solution of $f(x,y)=0$ at $O$.
Then by the assumption $I(B_3,B_3';O)=9$, we must have 
$\ord_x\, f'(x,\sum_{i=2}^\infty t_i x^i)=9.$
\begin{Lemma}
The sextic $C:=B_3\cup B_3'$ is of torus type if  and only if 
$t_2=0$. This implies that $O$ is a flex of $B_3$.
\end{Lemma}
{\em Proof.}
This assertion is given by Artal in \cite{Artal}. Our proof is
computational. In fact, if $t_2=0$, $O$ is a flex for both 
$B_3, B_3'$ and we see that 
$y=0$ is a flex tangent line  for $B_3,\,B_3'$.
Thus by Tokunaga's criterion, $y^2=0$ is the conic which gives a linear
torus decomposition. For the detail about linear torus decomposition,
we refer \cite{Reduced}.

Assume that $t_2\ne 0$ and we prove that any such $C$ is of
non-torus type. In fact, supposing $C$ to be a sextic of torus type, take
a torus decomposition
$f(x,y)f'(x,y)=f_2(x,y)^3+f_3(x,y)^2$.
Put $y_1:=y-\sum_{i=2}^\infty t_i x^i$.
So the assumption implies that 
\[
 f(x,y_1+\phi(x))\times f'(x,y_1+\phi(x))=y_1(y_1+cx^9)+\text{(higher
 terms)},\,c\ne 0
\]
and thus $y_1':=y-\sum_{i=2}^8\,t_i x^i$ is the maximal contact
coordinate
and it is also the solution of $f_3(x,y)=0$ in $y$ mod $x^9$ and 
$ f_2(x,\sum_{i=2}^8\,t_i x^i)\equiv 0\quad \modulo\, (x^6)$.
For the existence of a non-trivial conic $f_2(x,y)$, we see that
the coefficient must satisfy:
\begin{eqnarray}
J_0:=-3 t_4 t_3 t_2+2 t_3^3+t_5 t_2^2=0
\end{eqnarray}
(conics are five dimensional but we have 6 equation $coeff(
f_2(x,\sum_{i=2}^8\,t_i x^i),x,j)=0$
for $j=0,\dots, 5$).
Then we examine the other equations
\[
(\sharp)\quad \Coeff(g(x,\sum_{i=2}^8 t_i x^i),x,j)=0,\quad j=0,\dots, 8
\]
where $g(x,y)$ is a cubic which corresponds to either $f(x,y)$ or
$f'(x,y)$.
Write $g(x,y)$ as a generic cubic form with 10 coefficients (but by
scalar multiplication, one coefficient can be normalized to be 1, say
that of $ x^6$ is 1, and we have 9 free coefficients), we  solve the equations
($\sharp$) from $j=0$ to $j=0$ to $j=8$ to express the coefficients in
rational functions of $t_2,\dots, t_7$. At the last step,
we have one free coefficient undetermined, say the coefficient $c$ of $x^ay^b$ and a linear equation
$\Coeff(g(x,\sum_{i=2}^8 t_i x^i),x,8)=0$.
This is written as 
$K_1 c+K_0=0$
where $K_1,\, K_0$ are rational functions of $t_2,\dots, t_7$.
Thus to have two non-trivial cubic solutions, we need to have
$ K_1=K_0=0$.
Now we can easily check that there are no solutions
if we assume that $J_0(t_2,\dots, t_5 )=0$.
The other assertion will be checked
 in the explicit construction of examples.\qed

If $t_2\ne 0$, there is no line $L$ such that $L$ intersects only at
$O$.
Thus the sextic can not be of torus type by the classification in 
\cite {Reduced}. However the above argument is useful for the explicit
computation of non-torus sextics.

{\bf (1-a)} Let us consider the case $[A_{17}],\, [A_{17},A_1],\,
    [A_{17},2A_1]$.
The sextics of non-torus type with above configurations are obtained using 
above computation ($t_2\ne 0$).  Their Zariski partners are cubics intersecting at
    flex points.

(a-1) $[A_{17}],\, C=B_3+B_3'$, 
$C:=\{f_3(x,y)g_3(x,y)=0\}$ where
\begin{eqnarray*}
&{f_3} :=  - y^{3} - y^{2} + ( - x^{2} + x)\,y - x^{3} + x\\
&{g_3} := x - {\displaystyle \frac {10}{9}} \,y^{3} - x^{2}
\,y + {\displaystyle \frac {10}{9}} \,x\,y - y^{2} - 
{\displaystyle \frac {10}{9}} \,x^{3}
\end{eqnarray*}

(a-2)  $[A_{17}+A_1],\, C=B_3+B_3'$ with $B_3$ has a node:
$C:=\{f_3(x,y)g_3(x,y)=0\}$ where
\begin{eqnarray*}
&{f_3} := x - 3\,x\,y + 3\,x^{2}\,y - 2\,x^{2} - y^{2} + 2\,
y^{3} + 4\,y^{2}\,x + x^{3}\\
&{g_3} := x + x\,y + 11\,x^{2}\,y - 10\,x^{2} - y^{2} - 2\,y
^{3} + 8\,y^{2}\,x + 5\,x^{3}
\end{eqnarray*}

(a-3) $[A_{17}+2A_1],\, C=B_3+B_3'$ 
$C:=\{f_3(x,y)g_3(x,y)=0\}$, where cubics are nodal and
\begin{multline*}
{f_3} := {\displaystyle \frac {1}{48}} \,I(48\,x - 96\,x^{2
} - 264\,x\,y + 124\,y^{3} - 48\,y^{2} + 411\,y^{2}\,x + 264\,x^{
2}\,y + 48\,x^{3} \\
\mbox{} - 104\,I\,x\,y\,\sqrt{3} + 181\,I\,y^{2}\,x\,\sqrt{3} + 
104\,I\,x^{2}\,y\,\sqrt{3} + 16\,I\,x\,\sqrt{3} - 32\,I\,x^{2}\,
\sqrt{3} \\
\mbox{} + 16\,I\,x^{3}\,\sqrt{3} - 16\,I\,y^{2}\,\sqrt{3} + 52\,I
\,y^{3}\,\sqrt{3})\sqrt{3} \left/ {\vrule 
height0.41em width0em depth0.41em} \right. \!  \! ( - 1 + I\,
\sqrt{3}) \\
{g_3} :=  - {\displaystyle \frac {1}{8}} (48\,x - 72\,x^{2}
 + 10\,I\,x^{2}\,y\,\sqrt{3} - 25\,I\,y^{2}\,x\,\sqrt{3} + 56\,I
\,x\,y\,\sqrt{3} - 216\,x\,y + 68\,y^{3} - 48\,y^{2} \\
\mbox{} + 231\,y^{2}\,x + 138\,x^{2}\,y + 23\,x^{3} + 7\,I\,x^{3}
\,\sqrt{3} + 8\,I\,x^{2}\,\sqrt{3} - 16\,I\,x\,\sqrt{3} - 12\,I\,
y^{3}\,\sqrt{3} \\
\mbox{} + 16\,I\,y^{2}\,\sqrt{3}) \left/ {\vrule 
height0.41em width0em depth0.41em} \right. \!  \! ((3 + I\,\sqrt{
3})\,( - 1 + I\,\sqrt{3})) 
\end{multline*}

\vspace{.4cm}\noindent
\subsection{ Exceptional configuration: $[A_{11},2A_2,3A_1]$ with two cubic
components.}
In this case, we do  the similar computation. We compute sextics $C=B_3\cup
B_3'$
such that $B_3$ and $ B_3'$  have two $A_2$ at $(0,1)$ and $(1,0)$
respectively
and they intersect at $(0,-1)$ with intersection number $6$ to make
$A_{11}$. We have the following sextics 
of non-torus
type.
\begin{eqnarray*}
&f(x,y):=f_3(x,y)\, g_3(x,y)\\
&{f_3} :=  - {\displaystyle \frac {1}{44652}} (47\,\sqrt{3}
 + 168 + 195\,I + 74\,I\,\sqrt{3})(195\,y + 168\,I\,y^{3} - 168\,
I\,y^{2} - 168\,I\,y \\
&\mbox{} - 156\,I\,x^{2} + 74\,\sqrt{3}\,y^{3} - 74\,\sqrt{3}\,y^{
2} - 74\,\sqrt{3}\,y - 60\,\sqrt{3}\,x^{2} + 60\,\sqrt{3}\,y\,x^{
2} - 47\,I\,\sqrt{3}\,y^{3} \\
&\mbox{} + 47\,I\,\sqrt{3}\,y^{2} + 156\,I\,y\,x^{2} + 48\,I\,x^{2
}\,\sqrt{3} + 168\,x^{2} + 47\,I\,y\,\sqrt{3} - 48\,I\,y\,x^{2}\,
\sqrt{3} \\
&\mbox{} - 168\,y\,x^{2} + 195\,y^{2} - 195\,y^{3} - 195 + 168\,I
 + 74\,\sqrt{3} - 47\,I\,\sqrt{3} + 48\,x^{3}) 
\\
&g_3(x,y)=
{\displaystyle \frac {1}{35636460}} ( - 3033 - 1989\,I + 1313\,I
\,\sqrt{3} + 1361\,\sqrt{3})( - 90720\,x + 14898\,y - 15336\,y\,x
 \\
&\mbox{} + 3021\,\sqrt{3} - 52722\,x^{2} - 11749\,\sqrt{3}\,y^{3}
 - 8305\,\sqrt{3}\,y^{2} + 6465\,\sqrt{3}\,y + 3021\,\sqrt{3}\,x
^{2} \\
&\mbox{} - 1785\,\sqrt{3}\,y\,x^{2} + 438\,y\,x^{2} - 87573\,y^{2}
 - 24417\,y^{3} + 75384\,y^{2}\,x + 65388\,x^{3} \\
&\mbox{} - 6042\,x\,\sqrt{3} - 7839\,I\,y^{3} + 36333\,I\,y^{2} + 
60705\,I\,y + 16533\,I\,x^{2} + 5736\,I\,\sqrt{3} \\
&\mbox{} - 33066\,I\,x - 4680\,y\,x\,\sqrt{3} + 1362\,y^{2}\,x\,
\sqrt{3} - 6383\,I\,\sqrt{3}\,y^{3} + 22753\,I\,\sqrt{3}\,y^{2}
 \\
&\mbox{} + 231\,I\,y\,x^{2} + 5736\,I\,x^{2}\,\sqrt{3} + 34872\,I
\,y\,\sqrt{3} - 60936\,I\,y\,x - 11472\,I\,x\,\sqrt{3} \\
&\mbox{} - 27870\,I\,y^{2}\,x - 5532\,I\,y\,x^{2}\,\sqrt{3} - 
29340\,I\,y\,x\,\sqrt{3} - 17868\,I\,y^{2}\,x\,\sqrt{3} \\
&\mbox{} + 78054 + 16533\,I)
\end{eqnarray*}
\subsection{Examples of (b) and (d).} As the corresponding sextics of
torus type are linear, we only need to check the singularities are not 
colinear.

\noindent
{\bf (b) $\Si(C)\supset \{A_{11},A_5\}$.}
We put $A_{11}$  at (0,0) with tangent line $x=0$
and $A_5$  at (1,0).

\noindent
(b-1) $\Si(C)=[A_{11},A_5]:$
\begin{eqnarray*} 
f(x,y)=& (- y^{3} + (9\,x - 1)\,y^{2} + 7\,x^{3} - 8\,x^{2} + x)\times\\
&( - 2\,y^{3} + (5\,x - 1)\,y^{2} + ( - x^{2} + x)\,y + 4\,x^{3} - 
5\,x^{2} + x)
\end{eqnarray*}
(b-2) $\Si(C)=[A_{11},A_5,A_1]:$
\begin{multline*}
f(x,y)=
{\displaystyle \frac {1}{55}} ( - 16\,x\,y - 2\,y^{2} + 16\,x^{2}
\,y + 14\,y^{2}\,x - 8\,x^{2} + 5\,x^{3} + 3\,x) \\
(11\,x\,y + 4\,y^{2} - 11\,x^{2}\,y + 98\,y^{2}\,x - 5\,x^{2} + 
11\,x^{3} + 14\,y^{3} - 6\,x) 
\end{multline*}

\noindent
(b-3) $\Si(C)=[A_{11},A_5,2A_1]:$
\begin{multline*}
f(x,y) := ( - 175\,y^{3} + 11\,x^{2}\,\sqrt{3} - 88\,x^{2} - y^{3}\,
\sqrt{3} - 30\,y^{2}\,x\,\sqrt{3} - 18\,y\,x^{2}\,\sqrt{3} + 27\,
x + 61\,x^{3} + 48\,y^{2} \\
\mbox{} - 94\,y\,x + 83\,y\,x^{2} - 126\,y^{2}\,x + I\,x^{2}\,
\sqrt{3} + 5\,I\,y^{2}\,\sqrt{3} - 17\,I\,y^{2}\,x - I\,x\,\sqrt{
3} - 3\,I\,y\,x \\
\mbox{} + 11\,I\,\sqrt{3}\,y^{3} - 6\,y^{2}\,\sqrt{3} - 11\,x\,
\sqrt{3} - 8\,I\,x^{2} + 21\,I\,y^{2} + 8\,I\,x - 27\,I\,y^{3} + 
25\,I\,y^{2}\,x\,\sqrt{3} \\
\mbox{} + 8\,I\,y\,x\,\sqrt{3} + 15\,I\,y\,x^{2}\,\sqrt{3} + 2\,I
\,y\,x^{2} + 27\,y\,x\,\sqrt{3})(175\,y^{3} + 11\,x^{2}\,\sqrt{3}
 - 88\,x^{2} \\
\mbox{} + y^{3}\,\sqrt{3} - 30\,y^{2}\,x\,\sqrt{3} + 18\,y\,x^{2}
\,\sqrt{3} + 27\,x + 61\,x^{3} + 48\,y^{2} + 94\,y\,x - 83\,y\,x
^{2} \\
\mbox{} - 126\,y^{2}\,x - 3\,I\,y\,x + 11\,I\,\sqrt{3}\,y^{3} - 6
\,y^{2}\,\sqrt{3} - 11\,x\,\sqrt{3} - 27\,I\,y^{3} + 8\,I\,y\,x\,
\sqrt{3} \\
\mbox{} + 15\,I\,y\,x^{2}\,\sqrt{3} + 2\,I\,y\,x^{2} - 27\,y\,x\,
\sqrt{3} - 25\,I\,y^{2}\,x\,\sqrt{3} + 8\,I\,x^{2} - 21\,I\,y^{2}
 - 8\,I\,x \\
\mbox{} + I\,x\,\sqrt{3} - I\,x^{2}\,\sqrt{3} - 5\,I\,y^{2}\,
\sqrt{3} + 17\,I\,y^{2}\,x)
\end{multline*}

\vspace{.3cm}
{\bf (d) $\Si(C)\supset \{3A_5\}$, $C=B_3+B_3'$.}
We put three $A_5$ at $(0,1), \, (0,-1)$  and $(1,0)$.

\noindent
\begin{multline*}
[3A_5]:\,\mathit{f(x,y)} := ({\displaystyle \frac {3}{7}}  + x - y + y^{3
} - {\displaystyle \frac {3}{7}} \,x^{3} - x^{2} - 
{\displaystyle \frac {1}{7}} \,y\,x^{2} - y^{2}\,x - 
{\displaystyle \frac {3}{7}} \,y^{2}) \\
( - {\displaystyle \frac {4}{5}} \,x - y + 1 + {\displaystyle 
\frac {4}{5}} \,y^{2}\,x - {\displaystyle \frac {3}{5}} \,x^{2}
 + {\displaystyle \frac {3}{5}} \,y\,x^{2} + y^{3} - y^{2} + 
{\displaystyle \frac {2}{5}} \,x^{3}) 
\end{multline*}
\begin{multline*}
[3A_5,A_1]:\,\mathit{f(x,y)} := (y^{3} + {\displaystyle \frac {57}{4}} \,y^{2
}\,x + {\displaystyle \frac {1}{4}} \,y^{2} - y - y\,x + 
{\displaystyle \frac {1}{4}} \,x^{3} + {\displaystyle \frac {1}{4
}} \,x^{2} - {\displaystyle \frac {1}{4}} \,x - {\displaystyle 
\frac {1}{4}} ) \\
(y^{3} + {\displaystyle \frac {71}{5}} \,y^{2}\,x - y^{2} + 
{\displaystyle \frac {49}{5}} \,y\,x^{2} - 20\,y\,x - y + 
{\displaystyle \frac {67}{5}} \,x^{3} - {\displaystyle \frac {101
}{5}} \,x^{2} + {\displaystyle \frac {29}{5}} \,x + 1)
\end{multline*}
\begin{multline*}
[3A_5+2A_1]:\,\mathit{f(x,y)} :=  (21\,y^{2}\,x - 
12\,I\,\sqrt{3}\,y^{2}\,x + 12\,y\,x^{2} - I\,\sqrt{3}\,y\,x^{2}
 + 12\,y\,x - 2\,I\,\sqrt{3}\,y\,x - I\,y\,\sqrt{3} \\
\mbox{} + I\,y^{3}\,\sqrt{3} - 3 + 3\,y^{2} + 3\,x^{3} - 3\,x + 3
\,x^{2})( - 1 + 49\,x^{3} - 17\,y^{2}\,x + 84\,y\,x^{2} - 63\,x^{
2} \\
\mbox{} - 12\,y\,x + y^{2} + 15\,x - 7\,I\,\sqrt{3}\,y\,x^{2} - 6
\,I\,\sqrt{3}\,y\,x - 4\,I\,\sqrt{3}\,y^{2}\,x + I\,y\,\sqrt{3}
 - I\,y^{3}\,\sqrt{3})
\end{multline*}
\section{Three conics}
In this section, we study the last case $C=B_2+B_2'+B_2''$ with
the configuration of the singularities $[3A_5,3A_1]$.
Such a sextic is given when each pair of conics are intersecting at two
points: at one point, with intersection multiplicity 3 and at another
point, transversely. We can understand Zariski pairs
in this situation using conical flex
points.
Assume that the respective defining polynomials
of $B_2,\,B_2',\,B_2''$  are 
$f_2(x,y),\, g_2(x,y),\,h_2(x,y)$ and the location of two $A_5$'s are
$P_1=(0,1),\, P_2=(0,-1)$ with respective tangent cones are
$y\mp 1=0$.
We assume further $P_1\in B_2\cap B_2'$ and $P_2\in B_2\cap B_2''$.
We fix $B_2,\,B_2'$ generically and consider a linear system $\Phi$ of conics
$B_2''$ of dimension 2 such that $B_2''$ and  $B_2$
are tangent at $P_2$.
Under this situation we assert that
\begin{Proposition}
There exist 5 conical flex points $Q_i,\,i=1,\dots,5$ on $B_2'$
with respect to $\Phi$ so that $Q_1$ is a conical flex of torus type
and the other are of non-torus type.
\end{Proposition}
{\em Proof.}
To avoid the complexity of the equation, we choose a generic $B_2,B_2'$
so that 
\begin{eqnarray*}
&f_2(x,y)\,=\,(y^{2} - 1 + x^{2})\,\\
&g_2(x,y)\,=\, (y^{2} + ( - 
{\displaystyle \frac {2}{15}} \,\sqrt{130}\,x - {\displaystyle 
\frac {2}{3}} )\,y + {\displaystyle \frac {2}{3}} \,x^{2} + 
{\displaystyle \frac {2}{15}} \,\sqrt{130}\,x - {\displaystyle 
\frac {1}{3}} ) 
\end{eqnarray*}
We find 5 conical flex points on $B_2'$:
\begin{multline*}
Q_1=( - {\displaystyle \frac {1}{9}} \,\sqrt{130}, \,
{\displaystyle \frac {-17}{9}})\\
Q_2=(    {\displaystyle \frac {7}{3}} \,I\,\sqrt{10}\,
\sqrt{3} - {\displaystyle \frac {2}{3}} \,I\,\sqrt{13}\,\sqrt{10}
\,\sqrt{3},  - 4 + \sqrt{13} + {\displaystyle \frac {5}{3}} \,
I\,\sqrt{13}\,\sqrt{3} - {\displaystyle \frac {19}{3}} \,I\,
\sqrt{3})\\
Q3=( - {\displaystyle \frac {7}{3}} \,I\,\sqrt{10}\,
\sqrt{3} + {\displaystyle \frac {2}{3}} \,I\,\sqrt{13}\,\sqrt{10}
\,\sqrt{3}, - 4 + {\displaystyle \frac {19}{3}} \,I\,\sqrt{3}
 + \sqrt{13} - {\displaystyle \frac {5}{3}} \,I\,\sqrt{13}\,
\sqrt{3})\\
Q4=(- {\displaystyle \frac {7}{3}} \,I\,\sqrt{10}\,
\sqrt{3} - {\displaystyle \frac {2}{3}} \,I\,\sqrt{13}\,\sqrt{10}
\,\sqrt{3}, - 4 - {\displaystyle \frac {5}{3}} \,I\,\sqrt{13}
\,\sqrt{3} - {\displaystyle \frac {19}{3}} \,I\,\sqrt{3} - \sqrt{
13}),\,\\
Q5=( {\displaystyle \frac {7}{3}} \,I\,\sqrt{10}\,
\sqrt{3} + {\displaystyle \frac {2}{3}} \,I\,\sqrt{13}\,\sqrt{10}
\,\sqrt{3}, - 4 + {\displaystyle \frac {19}{3}} \,I\,\sqrt{3}
 - \sqrt{13} + {\displaystyle \frac {5}{3}} \,I\,\sqrt{13}\,
\sqrt{3})
\end{multline*}
$Q_1$ gives a sextic of torus type so that
$B_2''$ is given by
\begin{multline*}
h_2(x,y)\,=\, ({\displaystyle \frac {338}{201}} \,y - {\displaystyle \frac {104
}{1005}} \,y\,\sqrt{130}\,x - {\displaystyle \frac {104}{1005}} 
\,\sqrt{130}\,x + {\displaystyle \frac {137}{201}}  + y^{2} + 
{\displaystyle \frac {32}{201}} \,x^{2}) 
\end{multline*}
The other  conical flex points give sextics of non-torus type. For
example,
$Q_2$ gives $B_2''$ described as:
\begin{multline*}
h_2(x,y)=75 + 150\,I\,\sqrt{3} + 50\,\sqrt{13} - 45\,I\,
\sqrt{13}\,\sqrt{3} - 104\,\sqrt{13}\,\sqrt{10}\,x - 80\,\sqrt{10
}\,x \\
\mbox{} + 72\,I\,\sqrt{10}\,x\,\sqrt{3} + 300\,I\,y\,\sqrt{3} + 
790\,y - 90\,I\,y\,\sqrt{13}\,\sqrt{3} + 100\,y\,\sqrt{13} \\
\mbox{} - 104\,y\,\sqrt{13}\,\sqrt{10}\,x - 80\,x\,y\,\sqrt{10}
 + 72\,I\,x\,y\,\sqrt{10}\,\sqrt{3} - 45\,I\,y^{2}\,\sqrt{13}\,
\sqrt{3} + 50\,y^{2}\,\sqrt{13} \\
\mbox{} + 150\,I\,y^{2}\,\sqrt{3} + 715\,y^{2} + 320\,x^{2}
\end{multline*}


\begin{thebibliography}{10}

\bibitem{Artal}
E.~Artal~Bartolo.
\newblock Sur les couples des {Zariski}.
\newblock {\em J.\ Algebraic Geometry}, 3:223--247, 1994.

\bibitem{A-R-C-T}
E.~Artal~Bartolo, J.~Carmona, J.~I. Cogolludo, and H.-O. Tokunaga.
\newblock Sextics with singular points in special position.
\newblock {\em J. Knot Theory Ramifications}, 10(4):547--578, 2001.

\bibitem{ARCI}
E.~Artal~Bartolo, J.~Carmona~Ruber, and J.~I. Cogolludo~Agust{\'{\i}}n.
\newblock Essential coordinate components of characteristic varieties.
\newblock {\em Math. Proc. Cambridge Philos. Soc.}, 136(2):287--299, 2004.

\bibitem{AT2}
E.~Artal~Bartolo and H.-o. Tokunaga.
\newblock Zariski pairs of index 19 and {M}ordell-{W}eil groups of {$K3$}
  surfaces.
\newblock {\em Proc. London Math. Soc. (3)}, 80(1):127--144, 2000.

\bibitem{AT1}
E.~Artal~Bartolo and H.-o. Tokunaga.
\newblock Zariski {$k$}-plets of rational curve arrangements and dihedral
  covers.
\newblock {\em Topology Appl.}, 142(1-3):227--233, 2004.

\bibitem{NambaBook}
M.~Namba.
\newblock {\em Geometry of projective algebraic curves}.
\newblock Decker, New York, 1984.

\bibitem{ZariskiPairsI}
M.~Oka.
\newblock Zariski pairs on sextics 1.
\newblock {\em preprint, 2004}.

\bibitem{Okadual}
M.~Oka.
\newblock Geometry of cuspidal sextics and their dual curves.
\newblock In {\em Singularities---Sapporo 1998}, pages 245--277. Kinokuniya,
  Tokyo, 2000.

\bibitem{OkaAtlas}
M.~Oka.
\newblock Alexander polynomial of sextics.
\newblock {\em J. Knot Theory Ramifications}, 12(5):619--636, 2003.

\bibitem{Reduced}
M.~Oka.
\newblock Geometry of reduced sextics of torus type.
\newblock {\em Tokyo J. Math.}, 26(2):301--327, 2003.

\bibitem{Tokunaga-torus}
H.-o. Tokunaga.
\newblock (2,3) torus sextics and the {Albanese} images of 6-fold cyclic
  multiple planes.
\newblock {\em Kodai Math. J.}, 22(2):222--242, 1999.

\end{thebibliography}
\def\cprime{$'$} \def\cprime{$'$} \def\cprime{$'$} \def\cprime{$'$}
  \def\cprime{$'$}

\end{document}